\documentclass[graybox]{svmult}

\usepackage{mathptmx}       
\usepackage{helvet}         
\usepackage{courier}        
\usepackage{type1cm}        
\usepackage{makeidx}         
\usepackage{graphicx}        
\usepackage{difftrees}
\usepackage{amsmath,amssymb,amsfonts}
\usepackage{array}

\graphicspath{{contributed+books/templates/}} 
\usepackage{multicol}        
\usepackage[bottom]{footmisc}

\newcommand{\g}{\ensuremath{\mathfrak{g}}}
\newcommand{\e}{\ensuremath{\mathrm{e}}}
\newcommand{\Fr}{\ensuremath{\mathrm{Fr}}}
\newcommand{\R}{\mathbb{R}}
\renewcommand{\d}{\mathrm{d}}
\newcommand{\Id}{\mathrm{Id}}

\DeclareMathOperator{\ad}{ad}
\DeclareMathOperator{\Ad}{Ad}
\DeclareMathOperator{\dexp}{dexp}
\DeclareMathOperator{\grad}{grad}

\DeclareMathOperator{\SO}{SO}

\DeclareMathOperator{\GL}{GL}

\DeclareMathOperator{\St}{St}

\DeclareFontFamily{U}{MnSymbolC}{}
\DeclareSymbolFont{MnSyC}{U}{MnSymbolC}{m}{n}
\DeclareFontShape{U}{MnSymbolC}{m}{n}{
    <-6>  MnSymbolC5
   <6-7>  MnSymbolC6
   <7-8>  MnSymbolC7
   <8-9>  MnSymbolC8
   <9-10> MnSymbolC9
  <10-12> MnSymbolC10
  <12->   MnSymbolC12}{}
\DeclareMathSymbol{\intder}{\mathbin}{MnSyC}{'270}

\newcommand{\h}{\ensuremath{\mathfrak{h}}}
\newcommand{\Phis}{\overline{\Phi}}

\newcolumntype{x}{>{\centering\arraybackslash\hspace{0pt}}p{7mm}}

\let\amalg=\undefined
\let\coprod=\undefined
\DeclareSymbolFont{cmsymbols}{OMS}{cmsy}{m}{n}
\DeclareSymbolFont{cmlargesymbols}{OMX}{cmex}{m}{n}
\DeclareMathSymbol{\amalg}{\mathbin}{cmsymbols}{"71}
\DeclareMathSymbol{\coprod}{\mathop}{cmlargesymbols}{"60}

\makeindex             

\begin{document}

\title*{Lie group integrators}
\author{Brynjulf Owren}
\institute{Brynjulf Owren \at Norwegian University of Science and Technology, Alfred Getz vei 1, N-7491 Trondheim, Norway, \email{bryn@math.ntnu.no}
}
%
%
\maketitle

\abstract*{Each chapter should be preceded by an abstract (10--15 lines long) that summarizes the content. The abstract will appear \textit{online} at \url{www.SpringerLink.com} and be available with unrestricted access. This allows unregistered users to read the abstract as a teaser for the complete chapter. As a general rule the abstracts will not appear in the printed version of your book unless it is the style of your particular book or that of the series to which your book belongs.
Please use the 'starred' version of the new Springer \texttt{abstract} command for typesetting the text of the online abstracts (cf. source file of this chapter template \texttt{abstract}) and include them with the source files of your manuscript. Use the plain \texttt{abstract} command if the abstract is also to appear in the printed version of the book.}

\abstract{In this survey we  discuss a wide variety of aspects related to Lie group integrators. 
These numerical integration schemes for differential equations on manifolds have been studied in a general and systematic manner since the 1990s and the activity has since then branched out in several different subareas, focussing both on theoretical and practical issues.
From two alternative setups, using either frames or Lie group actions on a manifold, we  here introduce the most important classes of schemes used to integrate nonlinear ordinary differential equations on Lie groups and manifolds. We describe a number of different  applications where there is a natural action by a Lie group on a manifold such that our integrators can be implemented. An issue which is not well understood is the role of isotropy and how it affects the behaviour of the numerical methods. The order theory of numerical Lie group integrators has become an advanced subtopic in its own right, and here we give a brief introduction on a somewhat elementary level. Finally, we shall discuss Lie group integrators having the property that they preserve a symplectic structure or a first integral.
}

\section{Introduction}
\label{sec:1}
\cite{celledoni14ipo}.
Leonhard Euler is undoubtedly one of the most accomplished mathematicians of all times, and the modern theme 
\textit{Numerical methods for Lie groups} can be traced back to Euler in more than one sense.
Indeed, the simplest and possibly mostly used numerical approximation method for ordinary differential equations was first described in Euler's work {\em  Institutionum Calculi Integralis (1768), Volumen Primum, Ch VII)} and bears the name Euler's method.
And undoubtedly, the most used test case for Lie group integrators is the Euler's free rigid body system, which was derived in his amazing treatise on Mechanics in 1736.

In the literature on structural mechanics, Lie group integrators have been around for a long time, but  general, systematic studies of numerical integrators for differential equations on Lie groups and homogeneous manifolds began as recent as the 1990s. Some notable early contributions were those of Crouch and Grossman \cite{crouch93nio} and Lewis and Simo \cite{lewis94caf}. A series of papers by Munthe-Kaas \cite{munthe-kaas95lbt, munthe-kaas98rkm,munthe-kaas99hor} caused an increased activity from the late nineties when a large number of papers appeared over a short period of time. Many of these early results were summarised in a survey paper by Iserles et al. \cite{iserles00lgm}. The work on Lie group integrators has been inspired by many subfields of mathematics. Notably, the study of order conditions and backward error analysis uses results from algebraic combinatorics, Hopf algebras, and has more recently been connected to post Lie algebras by Munthe-Kaas and coauthors, see e.g. \cite{munthe-kaas13opl}. In connection with the search for inexpensive coordinate representations of Lie groups as well as their tangent maps, the classical theory of Lie algebras has been put to use in many different ways. The theory of free Lie algebras \cite{reutenauer93fla} has been used to find optimal truncations of commutator expansions for general Lie algebras, see  e.g. \cite{munthe-kaas99cia}. For coordinate maps taking advantage properties of a particular Lie algebra, tools such as root space decomposition \cite{owren00imb} and generalised polar decompositions \cite{krogstad09gpc} have been applied. Also there is of course a strong connection between numerical methods for Lie groups and the area of Geometric Mechanics. This connection is often used in the setup or formulation of differential equations in Lie groups or homogeneous manifolds where it provides a natural way of choosing a group action, and in order to construct Lie group integrators which are symplectic or conserve a particular first integral.

In this paper we shall discuss several aspects of Lie group integrators, we shall however not attempt to be complete. Important subjects related to Lie group integrators not covered here include the case of linear differential equations in Lie groups and the methods based on Magnus expansions, Fer expansions, and Zassenhaus splitting schemes. These are methods that could fit well into a survey on Lie group integrators, but for information on these topics we refer the reader to excellent expositions such as \cite{blanes09tme,iserles11mea}.  Another topic we leave out here is that of stochastic Lie group integrators, see e.g. \cite{malham08slg}.
We shall also focus on the methods and the theory behind them rather than particular applications, of which there are many. The interested reader may check out the references \cite{bras13nao,bruls10otu,celledoni14ait,kobilarov09lgi,terze15lgi,wensch09mis}.

In the next section we shall define a compact setup of notation for differential equations on differential manifolds with a Lie group action. Then in Section~\ref{sec:typesofschemes} we discuss some of the most important classes of Lie group integrators and give a few examples of methods. Section~\ref{sec:choiceofactions} briefly treats a selection of group actions which are interesting in applications. Then in Section~\ref{sec:isotropy} we shall address the issue of isotropy in Lie group integrators, in particular how the freedom offered by the isotropy group can either be used to reduce the computational cost of the integrator or be used to improve the quality of the solution. We take our own look at order theory and expansions in terms of a generalised form of B-series in Section~\ref{sec:ordertheory}.
Finally, in sections \ref{sec:symplectic} and \ref{sec:intpres} we consider Lie group integrators which preserve a symplectic form or a first integral.

%


\section{The setup} \label{sec:setup}
Let $M$ be a differentiable manifold of dimension $d<\infty$ and let the set of smooth vector fields on $M$ be denoted $\mathcal{X}(M)$. Nearly everything we do in this paper is concerned with the approximation of the $h$-flow of a vector field $F\in\mathcal{X}(M)$ for some small parameter $h$ usually called the \emph{stepsize}. In other words, we approximate the solution to the differential equation
\begin{equation} \label{ODE}
      \dot{y} = \frac{\d}{\d{t}}y=\left.F\right|_y,\qquad F\in\mathcal{X}(M).
\end{equation}
 Crouch and Grossman \cite{crouch93nio} used a set of smooth frame vector fields $E_1,\ldots,E_\nu$, $\nu\geq d$ on $M$, assuming
$$
      \mathrm{span}(E_1|_x,E_2|_x,\ldots,E_\nu|_x)=T_xM,\quad\mathrm{for\ each}\ x\in M.
$$
It can be assumed that the frame vector fields are linearly independent as derivations of the ring $\mathcal{F}(M)$ of smooth functions on $M$, and we denote their $\mathbb{R}$-span by $V$.
Any smooth vector field $F\in\mathcal{X}(M)$  can be  represented by $\nu$ functions $f_i\in\mathcal{F}(M)$
\begin{equation} \label{oderepframes}
      F|_x = \sum_{i=1}^\nu f_i(x) E_i|_x
\end{equation}
where the $f_i$  are not necessarily  unique. We then have a natural affine connection
$$
      \nabla_F G = \sum_{i} F(g_i) E_i
$$
which is flat with constant torsion $\tau(\sum_j f_jE_j,\sum_i g_iE_i)=\sum_{i,j} f_j g_i [E_i,E_j]$. For later, we shall need the notion of a frozen vector field relative to the frame. The freeze operator $\Fr: M\times\mathcal{X}(M)\rightarrow V$ is defined as
$$
    \Fr (x,F):=F_x=\sum_i f_i(x) E_i
$$
We note that the torsion can be defined by freezing the vector fields and then take the Lie-Jacobi bracket, i.e.
\begin{equation*}
     \tau(F,G)|_x =[\Fr(x,F), \Fr (x,G)]
\end{equation*}

Another setup is obtained by using a Lie group $G$ acting transitively from the left on $M$ \cite{munthe-kaas99hor}. The Lie algebra of $G$ is denoted $\g$. Any vector field $F$ can now be represented via a map $f:M\rightarrow\g$ and the infinitesimal action
 $\rho:\g\rightarrow\mathcal{X}(M)$
\begin{equation} \label{F_and_rho}
      F|_x = \rho\circ f(x)|_x,\qquad \rho:\g\rightarrow\mathcal{X}(M),\quad 
      \rho(\xi)|_x=\left.\frac{\mathrm d}{\mathrm{d}t}\right|_{t=0}\exp(t\xi)\cdot x
\end{equation}
We note that the map $f$ is not necessarily unique.

\section{Types of schemes} \label{sec:typesofschemes}
There is now a large variety of numerical integration schemes available, typically formulated with either of the setups of the previous section. In what follows, we assume that a finite dimensional Lie group $G$ acts transitively on a manifold $M$ and the Lie algebra of $G$ is denoted $\g$.

\subsection{Schemes of Munthe-Kaas type} \label{subsec:munthe-kaas-type}
Using the second setup, a powerful way of deriving numerical integrators devised by Munthe-Kaas \cite{munthe-kaas99hor} is to
\begin{enumerate}
\item In a neighborhood $U\subset\g$ of $0$ introduce a local diffeomorphism $\psi: U\rightarrow G$, such that
$$
\psi(0)=1\in G,\qquad T_0\psi = \Id_\g
$$
 \item Observe that the map $\lambda_{y_0}(v)=\psi(v)\cdot y_0$ is surjective on a neighborhood of the initial value $y_0\in M$.
\item Compute the pullback of the vector field $F=\rho\circ f$ along $\lambda_{y_0}$
\item Apply one step of a standard numerical integrator to the resulting problem on $U$
\item Map the obtained approximation back to $M$ by $\lambda_{y_0}$
\end{enumerate}

Even though the idea is very simple, there are several difficulties that need to be resolved in order to obtain fast and accurate integration schemes from this procedure.

Observe that the derivative of $\psi$ can be trivialised by right multiplication of the Lie group, such that
$$
     T_u\psi = TR_{\psi(u)}\circ \mathrm{d}\psi_u,\quad\mathrm{d}\psi_u:\g\rightarrow\g.
$$
With this in mind we set out to characterise the vector field on $U\subset\g$, this is a simple generalisation of a result in \cite{munthe-kaas99hor}
\begin{lemma} \label{lemma:rkmk}
Let $M$ be a smooth manifold with a left Lie group action $\Lambda: G\times M\rightarrow M$ and let $\g$ be the Lie algebra of $G$. Let $\psi:\g\rightarrow G$ be a smooth map, $\psi(0)=1$.
 Fix a point $m\in M$, and set $\Lambda_m=\Lambda(\cdot,m)$ so that
$$
      \rho(\xi)|_m = T_1\Lambda_m(\xi)
$$
Suppose $F\in\mathcal{X}(M)$ is of the form
$$
    F|_m = \rho(\xi(m))|_m,\quad\mbox{for some}\ \xi:M\rightarrow\g.
$$
Define $\lambda_m(u)=\Lambda(\psi(u),m)$. Then there is an open set $U\subseteq\g$ containing $0$ such that the vector field $\eta\in\mathcal{X}(U)$ defined as
$$
     \eta|_u = \mathrm{d}\psi_u^{-1}(\xi\circ\lambda_m(u))
$$
is $\lambda_m$-related to $F$.
\end{lemma}

The original proof in \cite{munthe-kaas99hor} where $\psi=\exp$ was adapted to general coordinate maps in \cite{owren00imb}. One step of a Lie group integrator is obtained just by applying a classical integrator, such as a Runge--Kutta method, to the corresponding locally defined vector field $\eta$ on $\g$. A Runge--Kutta method with coefficients $(A,b)$ applied to the problem $\dot{y}=\eta(y)$ in a linear space is of the following form
\begin{align*}
    Y_r&=y_0+h\sum_{j=1}^s a_{rj} k_j,\ r=1,\ldots, s,\\
    k_r&= \left.\eta\right|_{Y_r},\ r=1,\ldots,s, \\
    y_1&= y_0 + h\sum_{r=1}^s b_r k_r
\end{align*}
Here $y_0$ is the initial value, $h$ is the  step size, and $y_1\approx y(t_0+h)$ is the approximate solution at time $t_0+h$. The parameter $s\geq 1$ is called the number of stages of the method. If the matrix $A=(a_{rj})$ is strictly lower triangular, then the method is called \emph{explicit}.
 The application of such a method to the transformed vector field of Lemma~\ref{lemma:rkmk} ca be written out as follows
\begin{align}
  u_r &= h\sum_{j=1}^s a_{ij} \tilde{k}_j, &
  k_r &= \xi\circ\lambda_{y_0} (u_r),  &
  \tilde{k}_r &= \mathrm{d}\psi_{u_r}^{-1}(k_r),\quad r=1,\ldots,s,\label{rkmkstages}\\
  v &= h\sum_{r=1}^s b_r \tilde{k_r} & y_1 &= \lambda_{y_0}(v)\label{rkmkstep}.
\end{align}
A major advantage of this approach compared to other types of Lie group schemes is its interpretation as a smooth change of variables which causes the convergence order to be (at least) preserved from the underlying classical integrator. As we shall see later, one generally needs to take into account additional order conditions to account for the fact that the phase space is not a linear space. 

\subsubsection{Choosing the exponential map as coordinates, $\psi=\exp$}

The first papers by Munthe-Kaas on Lie group integrators \cite{munthe-kaas95lbt,munthe-kaas98rkm,munthe-kaas99hor} all used $\psi=\exp$ as coordinates on the Lie group.
In this case, there are several difficulties that need to be addressed in order to obtain efficient implementations of the methods. One is the computation of the exponential map itself. For matrix groups, there are a large number of algorithms that can be applied, see e.g. \cite{moler78ndw,moler03ndw}. In \cite{celledoni00ate,celledoni01mft} the authors developed approximations to the exponential which exactly map matrix Lie subalgebras to their corresponding Lie subgroups. Another issue to be dealt with is the differential of the exponential map
$$
     \mathrm{d}\exp_u := T_{\exp(-u)}R_{\exp(-u)}\circ T_u\exp
$$

\begin{lemma} \label{lemma2}
(Tangent map of $\exp:\g\rightarrow G$). Let $u\in\g$, $v\in T_u\g\simeq\g$. Then
$$
    T_u\exp(v) = \left.\frac{\mathrm d}{\mathrm{d}s}\right|_{s=0}\exp(u+sv)=TR_{\exp(u)}\circ\dexp_u v
    =\dexp_u(v)\cdot\exp(u)
$$
where
$$
\dexp_u(v)=\int_0^1\exp(r\ad_u)(v)\,\mathrm{d}r=\left.\frac{e^z-1}{z}\right|_{z=\ad_u}(v)
$$
\end{lemma}

\smallskip

\begin{proof} Let    $y_s(t) = \exp(t(u+sv))$ such that 
$$
      T_u\exp(v) = \left.\frac{\mathrm d}{\mathrm{d}s}\right|_{s=0} y_s(1)
$$
But for now we differentiate with respect to $t$ to obtain
\begin{equation} \label{ysdot}
     \dot{y}_s :=  \frac{\mathrm{d}}{\mathrm{d}t} y_s(t) = (u+sv)y_s(t) 
\end{equation}
We also note that $y_s(t) = \exp(tu)+\mathcal{O}(s)$ as $s\rightarrow 0$.
From \eqref{ysdot} we then get
$$
\dot{y}-uy_s = sve^{tu} + \mathcal{O}(s^2)
$$
and the integrating factor $e^{-tu}$ yields
$$
     \frac{\mathrm d}{\mathrm{d}t}\left(e^{-tu}y_s\right) = se^{-tu}ve^{tu} + \mathcal{O}(s^2)
$$
Integrating, and using that $y_s(0)=\mathrm{Id}$, we get
$$
   y_s(t) = e^{tu} + s\int_0^t e^{ru}ve^{-ru}e^{tu}\,\mathrm{d}r + \mathcal{O}(s^2)
$$
and so
$$
  \left.\frac{\mathrm d}{\mathrm{d}s}\right|_{s=0} y_s(1)=\int_0^1e^{ru}ve^{-ru}\,\mathrm{d}r\,e^{u}
  =\int_0^1e^{r\,\ad_u}(v)\,\mathrm{d}r\;e^{u}
$$
where we have used the well-known identity $\Ad_{\exp(u)}=\exp(\ad_u)$
in the last equality. Formally, we can write 
$$
\int_0^1e^{r\,\ad_u}(v)\,\mathrm{d}r=
\int_0^1 e^{rz}|_{z=\ad_u}(v)\,\mathrm{d}r=\left.\frac{e^z-1}{z}\right|_{z=\ad_u}(v)
$$
\end{proof}

It is often useful to consider the $\dexp$-map as an infinite series of nested commutators
$$
\dexp_u(v) = (I+\frac1{2!}\ad_u+\frac{1}{3!}\ad_u^2+\cdots)(v)=
v+\frac12\,[u,v]+\frac16\,[u,[u,v]]+\cdots
$$
In \eqref{rkmkstages} it is the inverse of $\dexp$ which is needed.
Note that the function 
$$
    \phi_1(z) =   \frac{e^z-1}{z}
$$
is entire, this means that its reciprocal
$$
   \frac{z}{e^z-1}
$$
is analytic where $\phi_1(z)\neq 0$. In particular this means that $\frac{1}{\phi_1(z)}$ has a converging Taylor series about $z=0$ in the open disk $|z|< 2\pi$. This series expansion can be shown to be
$$
        \frac{z}{e^z-1}=1-\frac{z}{2}+\sum_{k=1}^\infty \frac{B_{2k}}{(2k)!} z^{2k}
$$
where $B_{2k}$ are the \emph{Bernoulli numbers}, the first few of them are: $B_2=\frac16$, $B_4=-\frac1{30}$, $B_6=\frac1{42}$, $B_8=-\frac{1}{30}$, $B_{10}=\frac5{66}$. The map
$$
v=\dexp_u^{-1}(w)\qquad (\text{whenever}\ w=\dexp_u(v))
$$
is given precisely as
\begin{equation}\label{dexpinv}
       \dexp_u^{-1}(w) = \left.\frac{z}{e^z-1}\right|_{z=\ad_u}(w)=w-\frac{1}{2}[u,w]
       +\frac{B_2}{2!}[u,[u,w]] +\cdots
\end{equation}
We observe from \eqref{rkmkstages} that one needs to compute $\dexp_{u_r}^{-1}(k_r)$ and
that each $u_r=\mathcal{O}(h)$. This means that one may approximate the series in \eqref{dexpinv} by a finite sum, 
$$
\mathtt{dexpinv}(u,w,m)=w-\frac12[u,w]+\sum_{k=1}^m \frac{B_{2k}}{(2k)!} \ad_u^{2k}(w).
$$
One has $\mathtt{dexpinv}(u,w,m)\in\g$ for every $m\geq 0$ and furthermore
$$
      \dexp_{u_r}^{-1}(k_r)-\mathtt{dexpinv}(u_r,k_r,m)=\mathcal{O}(h^{2m+1})
$$
As long as the classical integrator has order $p\leq 2m+1$,  the resulting Munthe-Kaas scheme will also have order $p$.
There exists however a clever way to substantially reduce the number of commutators that need to be computed in each step. Munthe--Kaas and Owren \cite{munthe-kaas99cia} realised that one could form linear combinations of the stage derivatives $\tilde{k}_r$ in \eqref{rkmkstages} such that
$$
      Q_r = \sum_{j=1}^r \sigma_{r,j} \tilde{k}_j = \mathcal{O}(h^{q_r})
$$
for $q_r$ as large as possible for each $r$. Then these new quantities $Q_r$ were each given the grade $q_r$ and one considered the graded free Lie algebra based on this set. The result was a significant reduction in the number of commutators needed. Also Casas and Owren \cite{casas03cel}
provided a way to organise the commutator calculations to reduce even further the computational cost. Here is a Runge--Kutta Munthe--Kaas method of order four with four stages and minimal set of commutators
\begin{align*}
   k_1 &= h f(y_0),\\
   k_2 &= h f(\exp(\tfrac{1}{2}k_1) \cdot y_0), \\
   k_3 &= h f(\exp(\tfrac{1}{2}k_2-\tfrac{1}{8}[k_1,k_2])\cdot y_0), \\
   k_4 &= h f(\exp(k_3)\cdot y_0), &  \\
   y_1 &= \exp(\tfrac{1}{6}(k_1+2k_2+2k_3+k_4-\tfrac12[k_1,k_4]))\cdot y_0.
\end{align*}
For later reference, we also give the Lie-Euler method, a first order Lie group integrator generalising the classical Euler method
\begin{equation}\label{eq:Lie--Euler}
      y_1 = \exp(hf(y_0))\cdot y_0
\end{equation}

\subsubsection{Canonical coordinates of the second kind}

The exponential map is generally expensive to compute exactly. For matrix Lie algebras
$\g\subseteq\mathfrak{gl}(n,\mathbb{F})$ where $\mathbb{F}$ is either $\mathbb{R}$ or $\mathbb{C}$, standard software for computing $\exp$ numerically has a computational cost of $n^3$ to the leading order, and the constant in front of $n^3$ may be as large as $20-30$.  Another, yet completely general alternative to the exponential function is constructed as follows: Fix a basis for $\g$, say
$e_1,\ldots,e_d$ and consider the map
\begin{equation}\label{cc2}
\psi: v_1e_1+\cdots+v_de_d\mapsto\exp(v_1e_1)\cdot\exp(v_2e_2)\cdots\exp(v_de_d)
\end{equation}
Although it might seem unnatural to replace one exponential by many, one needs to keep in mind that if the basis can be chosen such that its exponential can be computed explicitly, it may still be an efficient method. For instance, in the general linear matrix Lie algebra $\mathfrak{gl}(m\mathbb{F})$ one may use the basis to be $e_{ij}=\mathbf{e}_i\mathbf{e}_j^T$ where $\mathbf{e}_i$ is the $i$th canonical unit vector in $\mathbb{R}^n$. Then
$$
\exp(\alpha e_{ij})=1 + \alpha e_{ij},\ i\neq j,\qquad
\exp(\alpha e_{ii})=1 + (\e^{\alpha}-1) e_{ii}
$$
So computing \eqref{cc2} takes approximately $n d$ operations which is much cheaper than computing the  exponential of a general matrix.

The difficulty lies however in computing the map $\mathrm{d}\psi_u^{-1}$ in an efficient manner. 
A method for this was developed in \cite{owren00imb}. The methodology is slightly different for solvable and semisimple Lie algebras. We here outline the main idea, for details we refer to 
the original paper \cite{owren00imb}. Differentiate \eqref{cc2} to obtain
$$
\mathrm{d}\psi_u(v) = v_1 e_1 + \sum_{i=2}^d\,v_i\;\Ad_{\e^{u_1e_1}}\circ\cdots\circ\Ad_{\e^{u_{i-1}e_{i-1}}}(e_i)
$$
The main idea is to find an equivalent expression which is a composition of cheaply invertible operators.
For this, we introduce a projector onto the span of the last $d-k$ basis vectors as follows
$$
     P_k: \sum_{i=1}^d v_ie_i\mapsto \sum_{i=k+1}^d v_ie_i
$$
where we let $P_0$ and $P_d$ equal the identity operator and zero operator on $\g$ respectively.
We may now define a modified version of the $\Ad$-operator, 
for any $u=\sum u_ie_i\in\g$, let
$$
\widehat{\Ad}_{\e^{u_ke_k}}=(\Id-P_k)+\Ad_{\e^{u_ke_k}}P_k
$$
This is a linear operator which acts as the identity operator on basis vectors $e_i,\ i\leq k$, and on basis vectors $e_i,\ i\geq k$ it coincides with  $\Ad_{\e^{u_ke_k}}$

\begin{definition}
An ordered basis $(e_1,\ldots,e_d)$ is called and admissible ordered basis (AOB) if, for each
$u=\sum u_je_j\in\g$ and for each $i=1,\ldots,d-1$, we have
\begin{equation}\label{aob}
\Ad_{\e^{u_1e_1}}\circ\cdots\circ\Ad_{\e^{u_ie_i}}P_i=
\widehat{\Ad}_{\e^{u_1e_1}}\circ\cdots\circ\widehat{\Ad}_{\e^{u_ie_i}}P_i
\end{equation}
\end{definition}
This definition is exactly what is needed to write $\mathrm{d}\psi_u$ as a composition of operators
\begin{proposition}
If the basis $(e_1,\ldots,e_d)$ is an AOB, then
$$
    d\psi_u = \widehat{\Ad}_{\e^{u_1e_1}}\circ\cdots  \widehat{\Ad}_{\e^{u_de_d}}
$$
\end{proposition}
Another important simplification can be obtained if an abelian subalgebra $\h$ of dimension $d-d^*$ can be identified. In this case the ordered basis can be chosen such that
$\h=\mbox{span}(e_{d^*+1},\ldots,e_d)$. Then $\left.\Ad_{\e^{u_ie_i}}\right|_{\h}$ for $i>d^*$ is the identity operator and therefore $\widehat{\Ad}_{\e^{u_ie_i}}$ is the identity operator on all of $\g$.
Summarizing, we have the following expression
$$
\mathrm{d}\psi_u^{-1} = \widehat{\Ad}^{-1}_{\e^{u_{d^*}e_{d^*}}}\circ\cdots\circ
 \widehat{\Ad}^{-1}_{\e^{u_{1}e_{1}}}
$$
Choosing typically a basis consisting of nilpotent elements, the inversion of each
$\widehat{\Ad}_{\e^{u_ie_i}}$ can be done cheaply by making use of the formula
$$
   \Ad_{\e^{u_ie_i}}=1+\sum_{k=1}^K \frac{u_i^k}{k!}\ad_{e_k},\quad\ad_{e_i}^{K+1}=0.
$$
For choosing the basis one may, for semisimple Lie algebras, use a basis known as the Chevalley basis. This arises from the root space decomposition of the Lie algebra
\begin{equation} \label{rootspacedecomposition}
\g=\h\oplus\coprod_{\alpha\in\Phi} \g_{\alpha}.
\end{equation}
Here $\Phi$ is the set of roots, and $\g_{\alpha}$ is the one-dimensional subspace of $\g$ corresponding to the root $\alpha\in\h^*$, see e.g. Humphreys \cite{humphreys72itl}.
$\h$ is the maximal toral subalgebra of $\g$ and it is abelian. In the previous notation, the 
number of roots is $d^*$ and the dimension of $\h$ is $d-d^*$. The following result whose proof can be found in \cite{owren00imb} provides a tool for determining whether an ordered Chevalley basis is an AOB.
\begin{theorem}\label{theo:rtsaob}
Let $\{\beta_1,\ldots,\beta_{d_*}\},\ d_*=d-\ell$, be the set
of roots $\Phi$ for a semisimple Lie algebra $\g$.
Suppose that a Chevalley basis is ordered as
\[
  (e_{\beta_1},\ldots,e_{\beta_{d_*}},h_1,\ldots,h_\ell)
\]
where $e_{\beta_i}\in\g_{\beta_i}$, and $(h_1,\ldots,h_\ell)$ is a
basis for $\h$.
Such an ordered basis is an AOB if
\begin{equation}\label{rtsaob}
k\beta_i+\beta_s=\beta_m,\ m<i<s\leq d_*,
k\in\mathbb{N}\quad\Rightarrow\quad
\beta_m+\beta_n\not\in\Phis,\ m<n\leq i-1.
\end{equation}
Here $\Phis=\Phi\cup\{0\}$.
\end{theorem}

\begin{example}
As an example, we consider $A_\ell=\mathfrak{sl}(\ell+1,\mathbb{C})$, commonly realized as the set of $(\ell+1)\times(\ell+1)$-matrices with vanishing trace. The maximal toral subalgebra is then the set of diagonal matrices in $\mathfrak{sl}(\ell+1,\mathbb{C})$. The positive roots are denoted
$$
     \{\beta_{i,j},1\leq i\leq j\leq\ell\}.
$$
Letting $\mathbf{e}_i$ be the $i$th canonical unit vector in $\mathbb{C}^{\ell+1}$, the root space
  corresponding to $\beta_{i,j}$  has a basis vector
$$
     \mathbf{e}_i\mathbf{e}_{j+1}^T \in \g_{\beta_{i,j}},\quad 1\leq i\leq j\leq \ell
$$
whereas the negative roots are associated to the basis vectors
$$
\mathbf{e}_{j+1}\mathbf{e_i}^T \in \g_{-\beta_{i,j}}.
$$
As a basis for $\h$, one may choose the matrices $\mathbf{e}_i\mathbf{e}_i^T-\mathbf{e}_{i+1}\mathbf{e}_{i+1}^T,\ 1\leq i\leq\ell$. The remaining difficulty now is to choose an ordering of the basis so that an AOB results. As indicated earlier, the basis for $\h$ may be ordered as the last ones, i.e.
with indices ranging from $d^*+1=\ell^2+\ell+1$ to $d=\ell^2+2\ell$.
With the convention $e_{\beta}\in\g_{\beta},\ \beta\in\Phi$,
$\h=\mathrm{span}(e_{\h_1},\ldots,e_{\h_\ell})$, let
\[
B=(e_{\beta_{i_1,j_1}},\ldots,e_{\beta_{i_m,j_m}},
 e_{-\beta_{i_1,j_1}},\ldots,e_{-\beta_{i_m,j_m}},e_{\h_1},\ldots,e_{\h_\ell}),
\]
where $i_1\leq i_2\leq\cdots\leq i_m$ and  $m=\ell(\ell+1)/2$. One can then prove by using Theorem~\ref{theo:rtsaob} that $B$ is an AOB. 
\end{example}
Similar details as for $A_\ell$ were also given for the other classical Lie algebras, 
$B_\ell, C_\ell, D_\ell$ and the exceptional case $G_2$ in \cite{owren00imb}. Also the case of solvable Lie algebras was considered.

\subsubsection{Other coordinate maps and retractions}
We have discussed two ways to choose coordinates on a Lie group as an ingredient in the Lie group integrators, these are canonical coordinates of the first and second kind. These choices are generic in the sense that they can be used for any finite dimensional Lie group with a corresponding Lie algebra. But if one allows for maps $\Psi:\g\rightarrow G$ that may only work for particular Lie groups there might be more options. Considering subgroups of the general linear group, a common type are those that can be embedded in $GL(n,\mathbb{R})$ via quadratic constraints, i.e. matrix groups of the form
$$
        G=\{A\in GL(n,\mathbb{R}): A^TJA=J\},
$$
for some $n\times n$-matrix $J$. It $J=\Id$, the identity matrix, then $G=SO(n,\mathbb{R})$ whereas if 
$J$ equals the constant Poisson structure matrix, then we recover the symplectic group $SP(2d,\mathbb{R})$. 
The Lie algebra of such a group consists of matrices
$$
\g = \{ a\in \mathfrak{gl}(n,\mathbb{R}): a^T J + Ja = 0 \}.
$$
As an alternative to the exponential map, while still keeping a map of the form $A=\chi(a)$ where
$\chi(z)$ is analytic in a neighborhood of $z=0$ is the Cayley transformation
$$
      \chi(z) = \frac{1+z/2}{1-z/2}.
$$
In fact, for any function $\chi(z)$ such that $\chi(-z)\chi(z)=1$ one has
$$
     a^T J + Ja = 0\quad\Rightarrow \quad \chi(a)^TJ\chi(a)=J.
$$
General software for computing  $\chi(a)$ for an $n\times n$-matrix has a computational cost of $\mathcal{O}(n^3)$, but the constant in front of $n^3$ is much smaller than what is required for the exponential map. The computation of the (inverse) differential of the Cayley transformation is also relatively inexpensive to compute, the right trivialised version is
$$
   \mathrm{d}\chi_y(u) = (1-\frac{y}{2})^{-1} u (I+\frac{y}{2})^{-1},\qquad
     \mathrm{d}\chi_y^{-1}(v) =(I-\frac{y}{2})v(I+\frac{y}{2}).
$$

\paragraph{\bf Retractions} \label{retractions}
 In cases of Lie group integrators where the Lie group has much higher dimension than the manifold it acts upon the computational cost may become too high for doing arbitrary calculations in the Lie algebra which is the way the Munthe-Kaas methods work. An option is then to replace the Lie algebra by a {\em retraction} which is a map retracting the tangent bundle $TM$ of the manifold into $M$
$$
         \phi: TM\rightarrow M.
$$
We let $\phi_x$ be the restriction of $\phi$ to $T_xM$ and denote by $0_x$ the zero-vector in $T_xM$. Following
\cite{adler02nmo} we impose the following conditions on $\phi$
\begin{enumerate}
\item $\phi_x$ is smooth and defined in an open ball $B_{r_x}(0_x)\subset T_xM$ of radius $r_x$ around $0_x$
\item $\phi_x(v)=x$ if and only if $v=0_x$.
\item $T_{0_x}\phi_x=1_{T_xM}$.
\end{enumerate}
Thus, $\phi_x$ is a diffeomorphism from some neighborhood $\mathcal{U}$ of $0_x$ to its image $\mathcal{W}=\phi_x(\mathcal{U})\subset M$. 

I a similar way as for the Munthe-Kaas type methods, the idea is now to make a local change of coordinates, setting for a starting point $y_0$, 
$$
    y(t) = \phi_{y_0}(\sigma(t)),\quad,\sigma(0)=0.
$$
This implies
$$
    \dot{y} = T_{\sigma}\phi_{y_0}\dot{\sigma} = F\circ\phi_{y_0}(\sigma).
$$
In some neighborhood of $0_{y_0}\in T_{y_0}M$ we have
\begin{equation} \label{reteq}
      \dot{\sigma} = (T_\sigma\phi_{y_0})^{-1}F\circ\phi_{y_0} (\sigma).
\end{equation}
The ODE on the vector space $T_{y_0}M$ can be solved by a standard integrator, and the resulting approximation over one step $\sigma_1$ to \eqref{reteq} is mapped back to $y_1 = \phi_{y_0}(\sigma_1)$ and the succeeding step is taken in coordinates from the tangent space $T_{y_1}M$ etc.
This way of introducing local coordinates for computation is in principle very simple, though it does not take into account the representation of the vector field as is the case with the Munthe-Kaas and Crouch-Grossman frameworks. Several examples of computationally efficient retractions can be found in \cite{celledoni02aco}, for instance in the orthogonal group by means of (reduced) matrix factorisations. 


\paragraph{\bf Retractions on Riemannian manifolds}
Geodesics can be used to construct geodesics on a  Riemannian manifold.  We define
$$
     \phi_x(v) = \exp_{x}(c) = \gamma_v(1),
$$
where $\gamma_v(t)$ is the geodesic emanating from $x$ with $\dot{\gamma}(0)=v$. The map $\exp_x$ is defined and of maximal rank in a neighborhood of $0_x\in T_xM$. The derivative of $\phi_x$ is related to the Jacobi field satisfying the Jacobi equation, see e.g. \cite[p. 70-82]{chavel93rga}. Let $\nabla$ be the Levi--Civita connection with respect to the metric on $M$ and let
$\mathbf{R}$ be the curvature tensor. Consider the vector field $Y$ defined along the geodesic $\gamma$, $\gamma(0)=x,
\dot{\gamma}(0)=v$ satisfying the boundary value problem
$$
\nabla^2_tY + \mathbf{R}(\dot{\gamma},Y)\dot{\gamma} = 0,\quad Y(0)=0,\quad Y(1)=w.
$$
Then
$$
   (T_v\phi_x)^{-1}(w) = (\nabla_tY)(0).
$$
Of particular interest in many applications is the case when there is a natural embedding of the Riemannian manifold into Euclidean space, say $V=\mathbb{R}^n$. In this case, one has for every $x\in M$ a decomposition of $V=T_xM\oplus N_xM$, where
$N_xM$ is the orthogonal complement of $T_xM$ in $V$. We may define a retraction
$$
   \phi_x(v) = x + v + n_x(v),
$$
where  $n_x(v)$ is defined in such away that $\phi_x(v)\in M$ for every $v$ belonging to a sufficiently small neighborhood of 
$0_x\in T_xM$. The derivative can be computed as
$$
 W :=  T_v\phi_x(w) =  w + T_vn_x(w),\quad T_vn_x(w)\perp T_xM
$$
so the image of the derivative in $T_{\phi_x(v)}M$ is naturally split into components in $T_xM$ and $N_xM$.
Now we can just apply the orthogonal projector $\mathbf{P}_{T_xM}$ onto $T_xM$ on each side to obtain
$$
w = (T_v\phi_x)^{-1} W = \mathbf{P}_{T_xM}W
$$

\subsection{Integrators based on compositions of flows}
\subsubsection{Crouch-Grossman methods}
What characterises the Munthe-Kaas type schemes is that they can be represented via a change of variables. But there are also Lie group integrators which do not have this property.
Crouch and Grossman \cite{crouch93nio} suggested method formats generalising both Runge--Kutta methods and linear multistep methods which they expressed in terms of frame vector fields, $E_1,\ldots,E_\nu$ as follows. Using the notation of \cite{owren99rkm} we present here the Runge--Kutta version of Crouch--Grossman methods.
\begin{align}
Y_r&=\exp(h a_{r,s}F_s)\circ\cdots\circ\exp(ha_{r,1}F_1)\,y_0 \label{eqYr},\\
F_r &= \Fr(Y_r,F) = \sum_{i=1}^\nu f_i(Y_r)E_i \label{eqFr},\\
y_1 &= \exp(hb_sF_s)\circ\exp(hb_{s-1}F_{s-1})\circ\cdots\circ\exp(hb_1 F_1) y_0, \label{eqy1}
\end{align}
where we have assumed that the vector field on $M$ has been written in the form \eqref{oderepframes}.
Here, the method coefficients $a_{r,j}$ and $b_j$ correspond to the usual coefficients of Runge-Kutta methods. In fact, whenever the frame is chosen to be the standard basis of $\mathbb{R}^n$, the method reduces to the familiar Runge--Kutta methods. Note that regardless of the ordering of exponentials in (\ref{eqYr},\ref{eqFr}) the method will reduce to the same  method whenever the flows are commuting. But in general a reordering of flows will alter the behaviour from order $h^3$ on, this indicates that it is not sufficient to enforce only the standard order conditions for Runge-Kutta methods on the $a_{r,j}$ and $b_r$ coefficients.  For classical explicit Runge--Kutta methods it is possible to obtain methods of order $p$ with $s=p$ stages for $p=1,2,3,4$, but for order $p\geq 5$ it is necessary that $s>p$. For Crouch-Grossman methods one can obtain $p=s$ for $p=1,2,3$, but for $p=4$ it is necessary to have at least five stages.
Crouch and Grossman \cite{crouch93nio} devised Runge-Kutta generalisations of methods of orders up to three, and Owren and Marthinsen \cite{owren99rkm} gave also and example of an explicit method of order four. We here give a third order method with three stages
\begin{align*}
Y_1 &= y_0,   & F_1 &= \Fr(Y_1,F),\\
Y_2 &=\exp(\tfrac34 h F_1)\,y_0,  & F_2&=\Fr(Y_2,F),\\
Y_3 &=\exp(\tfrac{17}{108}h F_2)\exp(\tfrac{119}{216}hF_1)\,y_0,& F_3&=\Fr(Y_3,F), \\[1mm]
y_1 &=\exp(\tfrac{24}{17}hF_3)\exp(-\tfrac{2}{3}hF_2)\exp(\tfrac{13}{51}hF_1)\,y_0.
\end{align*}
\subsubsection{Commutator-free Lie group integrators}\label{subsubsec:commfree}
A disadvantage of the Crouch--Crossman methods is that they use a high number of exponentials or flow calculations which is usually among the most costly operations of the method.
In fact, for an explicit method with $s$ stages one needs to compute $(s+1)s/2$ exponentials.
To improve this situation, Celledoni et al. \cite{celledoni03cfl} proposed a generalisation of the Crouch--Grossman Runge--Kutta style method which they called commutator-free methods
\begin{align} \label{commfree1}
Y_r &= \exp(\sum_k \alpha_{r,J_r}^k hF_k)\cdots\exp(\sum_k \alpha_{r,1}^k F_k), & F_r&=\Fr(Y_r,F), \\[1mm]
y_1 &=\exp(\sum_k \beta_{J}^k hF_k)\cdots\exp(\sum_k \beta_{1}^k F_k). \label{commfree2}
\end{align}
The intention here was to choose the number of flow calculations as small as possible by minimising 
the $J_r$, $J$. One may also here conveniently define
\begin{equation} \label{classcoeff}
     a_r^k = \sum_{j=1}^{J_r} \alpha_{r,j}^k,\quad
     b^k = \sum_{j=1}^J \beta_{j}^k,
\end{equation}
which will be the corresponding classical Runge--Kutta coefficients used when $M$ is Euclidean space.
There are actually commutator-free methods of order four with four stages which need the calculation of only five exponentials in total per step. In fact, the fourth order method presented in \cite{celledoni03cfl} needs effectively only four exponentials, because it reuses one exponential from a previous stage.
Writing as before $F_r=\Fr(Y_r,F)$ the method reads
\begin{align}
 Y_1&= y_0, \nonumber \\ 
 Y_2&=\exp(\tfrac{1}{2}hF_1)\cdot y_0, \nonumber \\
 Y_3&=\exp(\tfrac{1}{2}hF_2)\cdot y_0 \nonumber \\
 Y_4&= \exp(hF_3-\tfrac{1}{2}hF_1)\cdot Y_2,\label{cfexample} \\
 y_{\frac{1}{2}} &= 
 \exp(\tfrac{1}{12}h(3F_1+2F_2+2F_3-F_4))\cdot y_0,  \nonumber\\
 y_{1} &= 
 \exp(\tfrac{1}{12}h(-F_1+2F_2+2F_3+3F_4))\cdot y_{\frac{1}{2}}. \nonumber
  \end{align}
Note in particular in this example how the expression for $Y_4$ involves
$Y_2$ and thereby one exponential calculation has been saved.
Details on order conditions for commutator-free schemes can be found in \cite{owren06ocf}.

\section{ Choice of Lie group actions} \label{sec:choiceofactions}
The choice of frames or Lie group action is not unique and may have a significant impact on the properties of the resulting Lie group integrator.
 It is not obvious which action is the best or suits the purpose in the problem at hand. Most examples we know from the literature are using matrix Lie groups~$G\subseteq \GL(n)$, but the choice of group action depends on the problem and the objectives of the simulation. We give here some examples of situations where Lie group integrators can be used. Some more details are given in
 \cite{celledoni14ait}.

\subsection{Lie group acting on itself by multiplication} In the case $M=G$, it is natural to use either left or right multiplication 
\[
    \Lambda_L(g,m) = g\cdot m\quad\text{or}\quad R_m(m) = m\cdot g^{-1},\qquad g,m\in G.
\]
are both left Lie group actions by $G$ on $G$.
The corresponding maps  $\rho_L$ and $\rho_R$ defined in \eqref{F_and_rho} are
$$
\left.\rho_L(\xi)\right|_m = T_1R_m(\xi) = \xi\cdot m,\quad
\left.\rho_R(\xi)\right|_m = -T_1L_m(\xi) = -m\cdot\xi
$$
For a vector field $F\in\mathcal{X}(G)$ we use functions $f$ or $\tilde{f}$, $M\rightarrow\g$.
Using matrix notation we have
 $F|_g = f(g)\cdot g$ or $F|_g = -g\cdot\tilde{f}(g)$. Note that $\tilde{f}(g)$ is related to 
$f(g)$ through the adjoint representation of $G$, $\Ad \colon G\rightarrow \operatorname{Aut}(\g)$,
 \[
      f(g) = -\Ad_g \tilde{f}(g),\qquad \Ad_g = T_1L_g\circ T_1R_g^{-1}.
\]

\subsection{The affine group and its use in semilinear PDE methods}
Consider the semilinear partial differential equation
\begin{equation} \label{eq:semilinear}
     u_t = Lu + N(u),
\end{equation}
where $L$ is a linear differential operator and $N(u)$ is some nonlinear map, typically containing derivatives of lower order than $L$. Discretising \eqref{eq:semilinear} we obtain a system of  $n_d$ ODEs, for some large $n_d$, $L$ becomes an $n_d\times n_d$-matrix, and $N \colon \R^{n_d}\rightarrow \R^{n_d}$ a nonlinear function. 
We follow \cite{munthe-kaas99hor} and introduce a Lie group action on $\R^{n_d}$ by some subgroup of the affine group represented as the semidirect product $G=\GL(n_d)\ltimes\R^{n_d}$.
The group product, identity $1_G$, and inverse are given as
\[
   (A_1,b_1)\cdot (A_2,b_2) = (A_1 A_2, A_1b_2+b_1),\quad
   1_G=(1_{GL},0),\quad (A,b)^{-1}=(A^{-1}, -A^{-1}b).
\]
where $1_{GL}$ is the $n_d\times n_d$ identity matrix, and $0\in\R^{n_d}$.
The action on $\R^{n_d}$ is
\[
     (A,b)\cdot x = Ax + b,\qquad (A,b)\in G,\ x\in\R^{n_d},
\]
and the Lie algebra and commutator are given as
\[
   \g = (\xi, c),\ \xi\in\mathfrak{gl}(n_d),\ c\in\R^{n_d},\quad
   [(\xi_1,c_1),(\xi_2,c_2)] = ([\xi_1,\xi_2], \xi_1c_2-\xi_2c_1+c_1).
\]
In many interesting PDEs, the operator $L$ is constant, so it makes sense
to consider the $n_d+1$-dimensional subalgebra $\g_L$ of $\g$ consisting of elements 
$(\alpha L, c)$ where $\alpha\in\R$, $c\in\R^{n_d}$, so that the
map $f \colon \R^{n_d}\rightarrow \g_L$ is given as
\[
      f(u) = (L,N(u)).
\]
One parameter subgroups are obtained through the exponential map as follows
\[
   \exp(t(L,b)) = (\exp(tL), \phi(tL) tb).
\]
Here the entire function $\phi(z)=(\exp(z)-1)/z$ familiar from the theory of exponential integrators appears. As an example, one could now consider the Lie--Euler method \eqref{eq:Lie--Euler} in this setting, which coincides with the exponential Euler method
\[
   u_{1} = \exp(h(L,N(u_0))\cdot u_0 = \exp(hL)u_0 + h\phi(hL)N(u_0).
\]
There is a large body of literature on exponential integrators, going almost half a century back in time, see~\cite{hochbruck10ei} and the references therein for an extensive account. A rather general framework for exponential integrators were defined and studied in terms of order conditions in \cite{berland05bao}.

\subsection{The coadjoint action and Lie--Poisson systems}
Lie group integrators for this interesting case were studied by Eng{\o} and Faltinsen \cite{engo01nio}.
Suppose $G$ is a Lie group and the manifold under consideration is the dual space~$\g^*$ of its Lie algebra~$\g$. The coadjoint action by $G$ on $\g^*$ is denoted 
$\Ad_g^*$ defined for any $g\in G$ as
\begin{equation} \label{eq:coadjoint-action}
\langle \Ad_g^*\mu, \xi\rangle = \langle \mu, \Ad_g\xi\rangle,\quad\forall\xi\in\g, \mu \in \g^{*},
\end{equation}
for a duality pairing $\langle {\cdot}, {\cdot} \rangle$ between $\g^*$ and $\g$.
It is well known (see e.g.\ section~13.4 in \cite{marsden99itm}) that mechanical systems formulated on the cotangent bundle $T^*G$
with a left or right invariant Hamiltonian can be reduced to a system
on $\g^*$ given as
\[
     \dot{\mu} = \pm\ad^*_{\frac{\partial H}{\partial\mu}} \mu,
\]
where the negative sign is used in case of right invariance. The solution to this system preserves coadjoint orbits, which makes it natural to suggest the group action
\[
    g\cdot\mu = \Ad_{g^{-1}}^*\mu,
\]
so that the resulting Lie group integrator also respects this invariant.
For Euler's equations for the free rigid body, the Hamiltonian is left invariant and the coadjoint orbits are spheres in $\g^*\cong\R^3$.

\subsection{Homogeneous spaces and the Stiefel and Grassmann manifolds} \label{subsec:homogeneous}
The situation when $G$ acts on itself by left of right multiplication is a special case of a homogeneous space \cite{munthe-kaas97nio}, where the assumption is only that $G$ acts transitively and continuously on some manifold $M$. Homogeneous spaces are isomorphic to the quotient~$G/G_x$ where $G_x$ is the
\emph{isotropy group} for the action at an arbitrarily chosen point $x\in M$
\[
       G_x = \{ h\in G \mid h\cdot x = x\}.
\]
Note that if $x$ and $z$ are two points on $M$, then by transitivity of the action, $z=g\cdot x$ for some $g\in G$. Therefore,
whenever $h\in G_z$ it follows that $g^{-1} \cdot h \cdot g\in G_x$ so isotropy groups are isomorphic by conjugation~\cite{bryant95ait}.
Therefore the choice of $x\in M$ is not important for the construction of the quotient. For a readable introduction to this type of construction, see \cite{bryant95ait}, in particular Lecture 3.

 A much encountered example is the hypersphere
$M=S^{d-1}$ corresponding to the left action by $G=\SO(d)$,
 the Lie group of orthogonal $d\times d$ matrices with unit determinant.
 One has $S^{d-1} = \SO(d)/\SO(d-1)$. 
 
 The Stiefel manifold $\St(d,k)$ can be represented by the set of $d\times k$-matrices with orthonormal columns. An action on this set is obtained by left multiplication by $G=\SO(d)$. Lie group integrators for Stiefel manifolds are extensively studied in the literature, see e.g.\ \cite{celledoni03oti,krogstad03alc}.
An important subclass of the homogeneous spaces is the symmetric spaces, also
obtained through a transitive action by a Lie group $G$, where $M=G/G_x$, but here one requires in addition that the isotropy subgroup is an open subgroup of the fixed point set of an involution of $G$~\cite{munthe-kaas01aos}.  
A prominent example of a symmetric space in applications is the Grassmann manifold, obtained as $\SO(d)/(\SO(k)\times \SO(d-k))$.

\subsection{Isospectral flows} In isospectral integration one considers dynamical systems evolving on the manifold of $d\times d$-matrices sharing the same Jordan form.
Considering the case of symmetric matrices, one can use the transitive group action by $\SO(d)$
given as
\[
       g\cdot m = g m g^T.
\]
This action is transitive, since any symmetric matrix can be diagonalised by an appropriately chosen orthogonal matrix.  If the eigenvalues are distinct, then the isotropy group is discrete and consists of all matrices in $\SO(d)$ which are diagonal.

Lie group integrators for isospectral flows have been extensively studied, see for example
\cite{calvo96rkm1,calvo97nso}.  See also~\cite{celledoni01ano} for an application to the KdV equation.

\subsection{Tangent and cotangent bundles}
For mechanical systems the natural phase space will often be the tangent bundle $TM$ as in the Lagrangian framework or the cotangent bundle $T^*M$ in the Hamiltonian framework.
The seminal paper by Lewis and Simo \cite{lewis94caf} discusses several Lie group integrators for mechanical systems on cotangent bundles, deriving methods which are symplectic, energy and momentum preserving.
Eng{\o}~\cite{engo03prk} suggested a way to generalise the Runge--Kutta--Munthe-Kaas methods into a partitioned version when $M$ is a Lie group. 
Marsden and collaborators have developed the theory of Lie group integrators from the variational viewpoint over the last two decades. See \cite{marsden01dma} for an overview.
For more recent work pertaining to Lie groups in particular, see
\cite{lee07lgv,bou-rabee09hpi,saccon09mrf}. 

\section{ Isotropy}\label{sec:isotropy}
In Section~\ref{sec:setup} it was pointed out that when using frames to express the vector field $F\in\mathcal{X}(M)$, the functions
$f_i: M\rightarrow M$ in the expression $F|_x=\sum_i f_i(x)E_i|_x$ were not necessarily unique.
Similarly, using a map  $f:M\rightarrow\g$ in the group action framework to write $F=\rho\circ f$ where $\rho$ is the infinitesimal action, we also remarked that $f$ is not generally unique. In fact, taking any other map $z: M\rightarrow\g$ satisfying
$\rho(z(x))|_x=0$, then $F|_x=\rho(f(x))|_x = \rho(f(x)+z(x))|_x$ hence $f+z$ represents the same vector field as $f$.
So even though adding the map $z$ to $f$ does not alter the vector field $F$, it \emph{does} generally alter any numerical Lie group integrator. A typical situation with isotropy was described in Section~\ref{subsec:homogeneous} with homogeneous manifolds.

An interesting example is the two-sphere $S^2\simeq SO(3)/SO(2)$. If we represent elements of $\mathfrak{so}(3)$ as vectors in 
$\mathbb{R}^3$ and points on the sphere as three-vectors of unit length, we may express \eqref{F_and_rho} by using the cross-product in $\mathbb{R}^3$
$$
     F|_y = f(y)\times y =(f(y)+\alpha(y)y)\times y
$$
where $\alpha:S^2\rightarrow\mathbb{R}$ is any smooth function. The Lie--Euler method \eqref{eq:Lie--Euler} with initial value $y(0)=y$ would read for the first step
$$
      y_1 = \exp(h(f(y)+\alpha(y)y))y
$$
We may  assume that $(f(y),y)=0$ for all $y\in S^2$ with the Euclidean inner product. Then, by using Rodrigues formula for the exponential and simplifying, we can write
$$
    y_1 = (1-h^2\frac{1-\cos\theta}{\theta^2})y+h\frac{\sin\theta}{\theta}f(y)\times y +h^2\frac{1-\cos\theta}{\theta^2}\alpha(y)f(y)
$$
where $\theta=h\sqrt{\|f(y)\|^2+|\alpha(y)|^2}$. For small values of the step size $h$ this can be approximated as
$$
     y_1\approx (1-\frac12 h^2)y + h f(y)\times y +\frac{1}{2}h^2\alpha(y) f(y)
$$
so we see that by choosing $\alpha=0$ one moves from $y$ in the direction of the vector field $f(y)\times y$ whereas any nonzero $\alpha$ will give a  contribution to the increment in the tangent plane orthogonal to the vector field.
Lewis and Olver \cite{lewis02gia} pursue this analysis much further and show how the isotropy parameter $\alpha(y)$ can be used to minimise the orbital error by requiring that the curvatures of the exact and numerical solution agree.

Another situation in which the special care should be taken in choosing the freedom left by isotropy is the case when $M$ is a homogeneous manifold and the acting group $G$ has much higher dimension than that of $M$ as may be the case for instance with Stiefel and Grassmann manifolds discussed in \ref{subsec:homogeneous}.  Typically, a naive implementation of a Lie group integrators has a computational cost of order $d^3$ whenever $G$ is represented as a subgroup of $GL(d,\mathbb{R})$. But it is then useful to choose the map $f: M\rightarrow\g\subset\mathfrak{gl}(d,\mathbb{R})$ in such a way that the resulting Lie algebra element has the lowest possible rank as a matrix. In \cite{celledoni03oti} this idea was applied to Stiefel manifolds $\St(d,k)\simeq SO(d)/SO(d-k)$ and it was shown that a clever choice of isotropy component in $f(y)$ results in $\text{rank} f(y)=2k$, and then it was shown that the complexity of the Lie group integrator could be reduced from $\mathcal{O}(d^3)$ to $\mathcal{O}(dk^2)$. Krogstad \cite{krogstad03alc} suggested a similar approach which also leads to $\mathcal{O}(dk^2)$ complexity algorithms.

\section{ Order theory for Lie group integrators} \label{sec:ordertheory}
Order theory is concerned with the convergence order in terms of the step size for numerical integrators.
In the setting of manifolds one can define the concept precisely as follows: Let $\Phi_{h,F}: M\rightarrow M$ be a numerical flow map applied to $F\in\mathcal{X}(M)$  in the sense that the numerical method is defined through $y_{n+1}=\Phi_{h,F}(y_n)$. We shall say that this method has \emph{order} $p$ if, for any $F\in\mathcal{X}(M)$, $\psi\in\mathcal{F}(M)$, and $y\in M$ it holds that
$$
     \psi(\exp(h F) y) - \psi(\Phi_h (y)) = \mathcal{O}(h^{p+1})
$$
where here and in what follows write simply $\Phi_h$ for $\Phi_{h,F}$.  The order theory amounts to writing the exact and numerical solutions to the problem in powers of the stepsize $h$ and comparing these term by term.  A systematic development of the order theory for Crouch-Grossman and commutator-free Lie group was done in \cite{owren99rkm} and \cite{owren06ocf} respectively, but the theory is now well established and founded on an advanced algebraic machinery  much thanks to Munthe-Kaas and  coauthors \cite{munthe-kaas08oth,lundervold11hao,munthe-kaas13opl,ebrahimi-fard15otl} and Murua \cite{murua99fsa,murua06tha}.
 We skip most of the details here, and remind the reader that what we refer to as "Munthe-Kaas like" schemes in~\ref{subsec:munthe-kaas-type} are relatively easy to deal with in practice since they can be seen as the application of a standard "vector space" Runge--Kutta schemes under a local change of variables. The schemes based on compositions of flows are less straightforward, one here needs to use the theory of ordered rooted trees. A key formula is the series giving the pullback of an arbitrary  function $\psi\in\mathcal{F}(M)$ along the flow of a vector field
$$
\exp(hF)^*\psi = \psi\circ\exp(hF)= \psi + hF(\psi) + \frac{h^2}{2!} F^2(\psi) + \cdots
$$
Note also that by Leibniz' rule
$$
     F^2 = \Fr(\cdot,F)^2 + \nabla_F(F) = \sum_{i,j}( f_j f_i E_jE_i + f_jE_j(f_i)E_i)
$$
The two terms can be associated to ordered rooted trees, each  having  three nodes.
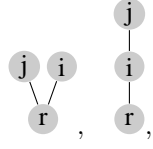
\tikzstyle mytree=[level distance=7mm,
every node/.style={scale=1.0,fill=gray!40,circle,minimum size=4mm,inner sep=1pt},
level 1/.style={sibling distance=5mm,nodes={fill=gray!40}},
level 2/.style={sibling distance=5mm,nodes={fill=gray!40}},
level 3/.style={sibling distance=5mm,nodes={fill=gray!40}}]
$$
\begin{tikzpicture}[mytree]
\node  {r} [grow=up]
child {node {i}}
child {node{j}};
\end{tikzpicture},\quad
\begin{tikzpicture}[mytree]
\node{r} [grow=up]
child {node{i}
child {node{j}}};
\end{tikzpicture},
$$
The reader may find it unnatural to include the fictitious root node $r$, the alternative is to represent the terms by means of ordered forests of rooted trees, but we shall stick with this notation, used for instance by Grossman and Larson \cite{grossman05das}, since it simplifies parts of the presentation.
 
Similarly, one can compute $F^3$
\begin{align}\label{F-cubed}
F^3 &= \sum_{i,j,k}( f_kf_jf_i E_kE_jE_i + f_kf_jE_k(f_i)E_jE_i + f_kE_k(f_j)f_i E_jE_i\\& \nonumber
+f_kf_jE_j(f_i)E_kE_i+f_kf_jE_kE_j(f_i)E_i+f_kE_k(f_j)E_j(f_i)E_i
\end{align}
where each term is associated to ordered trees as follows
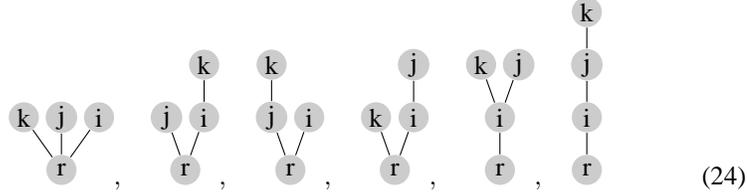
\begin{equation} \label{F-cubed-trees}
\begin{tikzpicture}[mytree]
\node  {r} [grow=up]
child {node {i}}
child {node{j}}
child {node{k}};
\end{tikzpicture},\quad
\begin{tikzpicture}[mytree]
\node {r} [grow=up]
child {node {i} child {node {k}}}
child {node {j}};
\end{tikzpicture},\quad
\begin{tikzpicture}[mytree]
\node {r} [grow=up]
child {node {i}}
child {node {j} child {node {k}}};
\end{tikzpicture},\quad
\begin{tikzpicture}[mytree]
\node {r} [grow=up]
child {node {i} child {node {j}}}
child {node {k}};
\end{tikzpicture},\quad
\begin{tikzpicture}[mytree]
\node {r} [grow=up]
child {node {i}
child {node{j}}
child {node{k}}};
\end{tikzpicture},\quad
\begin{tikzpicture}[mytree]
\node{r} [grow=up]
child {node{i}
child {node{j}
child {node{k}}}};
\end{tikzpicture}
\end{equation}
Going from $F^2$ to $F^{3}$ amounts to the following in our drawing convention: For every tree in $F^2$ run through all of its nodes and add a new node with a new label ($k$) as a new 'leftmost' child.
 Assuming the labels are ordered, say $i<j<k$, each labelled, monotonically ordered tree appear exactly once. In the convention used above we the children of any node will have increasing labels from right to left. Also they must necessarily increase from parent to child. With this convention, it is transparent that the only difference between the second and the fourth term is that the summation indices $j$ and $k$ have been interchanged so these two terms are the same. In general, we have an equivalence relation on the set of labelled ordered trees. Let the set of $q+1$ labelled nodes with a total order on the labels be denoted $\mathcal{A}_{q+1}$. To an ordered tree with $q+1$ nodes, we associate a child-to-parent map 
$\mathbf{t}:\mathcal{A}_{q+1}\backslash\{r\}\rightarrow\mathcal{A}_{q+1}$.
Now we define two labelled forests with child-to-parent maps $\mathbf{u}$ and $\mathbf{v}$ to be equivalent if there exists a permutation $\sigma:\mathcal{A}_{q+1}\rightarrow\mathcal{A}_{q+1}$ such that
\begin{enumerate}
\item $\sigma(r)=r$
\item $\mathbf{u}\circ\sigma = \sigma\circ\mathbf{v}$ on $\mathcal{A}_{q+1}\backslash\{r\}$
\item For all $z_1, z_2\in\mathcal{A}_{q+1}\backslash\{r\}$ such that $\mathbf{u}(z_1)=\mathbf{u}(z_2)$ one has
$$
   \sigma(z_1) < \sigma(z_2)\Rightarrow z_1<z_2
$$
\end{enumerate}
The ordered rooted trees are from now on taken to be equivalence classes of their labelled counterparts.
\begin{svgraybox}
 \textbf{Notation summary.} 
 \begin{description}
 \item[$T_O$]  is the set of ordered rooted trees
 \item[$F_O$] is the set of ordered forests, i.e. an ordered set of trees, $t_1t_2\ldots t_\mu$, each $t_i\in T_O$
 \item[$|t|$] is the total number of nodes in a tree (resp. forest)
  \item[$\bar{T}_O$] is the set $T_O$ augmented by the 'empty tree', $\emptyset$ where $|\emptyset|=0$, similarly $\bar{F}_O=F_O\cup\emptyset$, i.e. $F_O$ augmented with the empty forest.
 \item[$B_+$] is a map that takes any ordered forest into a tree by joining the roots of each tree to a new common node,
     $|B_+(t_1\ldots t_\mu)| = 1+|t_1\ldots t_\mu|$, we also let $B_+(\emptyset)=\ab$
\item[$B_-$] is a map $B_-: T_O\rightarrow \bar{F}_O$ that removes the root of a tree and leaves the forest of subtrees.
$B_+\circ B_-$ is the identity map on $T_O$.
\item[$\mathcal{T}_O$] The $\mathbb{R}$-linear space having the elements of $T_O$ as basis.
 \end{description}
\end{svgraybox}

\noindent We may now define $\alpha(t)$ to be the number of elements in the equivalence class which contains the labelled ordered tree $t$.
A formula for $\alpha(t)$ for ordered rooted trees was found in \cite{owren99rkm}.
\begin{proposition} \label{prop:alpha}
Set $\alpha(\ab)=1$. For any $t=B_+(t_1t_2\ldots t_\mu)\in T_O$ 
$$
\alpha(t) = \prod_{\ell=1}^\mu
\left(
\begin{array}{cc}
\sum_{i=1}^{\ell} |t_i|-1\\
|t_{\ell}|-1
\end{array}
\right) \alpha(t_\ell)
$$
where $|t|$ is the total number of nodes in the tree $t$. 
\end{proposition}
\begin{example}
We can work out the coefficient for a couple of examples. For instance
$$
\alpha(\aabb)=\alpha(B_+(\ab)) = \left(\begin{array}{c} 0\\0\end{array}\right)\alpha(\ab)=1
$$
For the tree $\aabaabbb$ we get
$$
\alpha(\aabaabbb) = \left(\begin{array}{c} 0\\0\end{array}\right)  \left(\begin{array}{c} 2\\1\end{array}\right)
\alpha(\ab)\alpha(\aabb)= 1\cdot 2\cdot 1\cdot 1 = 2,
$$
whereas
$$
\alpha(\aaabbabb) =  \left(\begin{array}{c} 1\\1\end{array}\right)  \left(\begin{array}{c} 2\\0\end{array}\right)
\alpha(\aabb)\alpha(\ab) = 1\cdot 1\cdot 1\cdot 1 = 1.
$$
\end{example}
For any tree $t=B_+(t_1\ldots t_\mu)$, $|t|>1$, and vector field $F$ we can associate operators $\mathbf{F}(t)$ acting on functions of $M$, defined in a recursive manner, setting
\begin{align*}
   \mathbf{F}(\ab)&=1,\\
   \mathbf{F}(t)&=\sum_{i_1\ldots i_\mu} \mathbf{F}(t_1)(f_{i_1})\cdots\mathbf{F}(t_\mu)(f_{i_\mu})
   E_{i_1}\ldots E_{i_\mu}
\end{align*}
Note in particular that according to the definition
$$ 
     \mathbf{F}(\aabb) = \sum_{i} f_i E_i = F.
$$
Similarly, there is a counterpart to the frozen vector fields relative to a point $x\in M$
$$
     \mathbf{F}_x(t) = \sum_{i_1\ldots i_\mu} \mathbf{F}(t_1)(f_{i_1})|_x\cdots\mathbf{F}(t_\mu)(f_{i_\mu})|_x
      E_{i_1}\ldots E_{i_\mu}
$$

The formal Taylor expansion of the flow of the vector field $F$ is 
$$
\psi\circ \exp(hF)y = \sum_{k=0}^\infty \frac{h^k}{k!} \left.F^k(\psi)\right|_y,\quad \psi\in\mathcal{F}(M)
$$
and in view of the above discussion we can therefore write the formal expansion
\begin{equation} \label{Bseries:exact}
     \psi\circ \exp(hF)y = \sum_{t\in T_O} \frac{h^{|t|-1}}{(|t|-1)!}\,\alpha(t)\,\mathbf{F}(t)(\psi)|_y
\end{equation}
This infinite series is in fact a special instance of what we may call a Lie-Butcher series, in which the coefficient $\alpha(t)/|t|!$ is replaced by a general map $\mathbf{a}:F_O\rightarrow\mathbb{R}$. We define the operator series
\begin{equation} \label{Bseries}
     B(\mathbf{a},x)=\sum_{t\in T_O} h^{|t|-1}\mathbf{a}(t)\mathbf{F}_x(t)
\end{equation}

At this point it is convenient to consider the free $\mathbb{R}$-vector space $\mathcal{T}_O$ over the set of  ordered rooted trees $T_O$ and to extend the map $\mathbf{F}$ to a linear map between $\mathcal{T}_O$ and  the space of higher order derivations indicated above.
The algebraic structures on $\mathcal{T}_O$ by now well-known from the literature, see e.g. \cite{berland05aso,lundervold11hao,munthe-kaas13opl} stem from the algebra on higher order derivations under composition by requiring $\mathbf{F}$ to be a algebra homomorphism. We will not further pursue these issues here, but instead we shall briefly consider how the commutator-free methods of Section~\ref{subsubsec:commfree} can be expanded in a series with the same type of terms as the exact flow, the details can be found in \cite{owren06ocf}. We introduce the concatenation product on $\mathcal{T}_O$ which is linear in both factors and on single trees $u$ and $v$ defined as 
\begin{equation} \label{concat}
u\cdot v =B_+( u_1u_2\ldots u_\mu v_1\ldots v_\nu),\quad u= B_+(u_1u_2\ldots u_\mu),\  v=B_+(v_1v_2\ldots v_\nu)
\end{equation}
This is a morphism under $\mathbf{F}_x$ since obviously $\mathbf{F}_x(u)\circ \mathbf{F}_x(v)=\mathbf{F}_x(u\cdot v)$.
\subsection{Order conditions for commutator-free Lie group integrators}
We now provide a few simple tools which we formulate through some lemmas with easy proofs which we omit, consult \cite{owren06ocf} for details
\begin{lemma}
Suppose $\phi_{\mathbf{a}}$ and $\phi_{\mathbf{b}}$ are maps of $M$ with B-series $B(\mathbf{a},\cdot)$ and $B(\mathbf{b},\cdot)$ respectively, such that $\mathbf{a}(\ab)=\mathbf{b}(\ab)=1$. This means that formally, for any smooth function $\psi$
$$
\psi(\phi_{\mathbf{a}}(y))=B(\mathbf{a},y)(\psi)|_y,\qquad \psi(\phi_{\mathbf{b}}(y))=B(\mathbf{b},y)(\psi)|_y
$$
The the composition of maps $\phi_a\circ\phi_b$ has a B-series $B(\mathbf{ab},y)$ with coefficients
$$
   \mathbf{ab}(\ab)=1,
$$
and for $t=t_1t_2\ldots t_\mu$
$$
\mathbf{ab}(t)=\sum_{u\cdot v=t} \mathbf{a}(v)\mathbf{b}(u)=\sum_{k=0}^\mu \mathbf{a}(t_{k+1}\ldots t_\mu)
\mathbf{b}(t_1\ldots t_k).
$$
\end{lemma}
\begin{lemma}\label{lemma:BseriesFrozen}
Suppose $a=\phi_{\mathbf{a}}(x)$ has a B-series $B(\mathbf{a},x)$. Then the frozen vector field
 $F_a=\sum_i f_i(a)E_i\in V$ has the B-series $hF_a=B(\mathbf{F}_a,x)$ where
 \begin{align*}
 \mathbf{F}_a(\ab)&=0 \\
 \mathbf{F}_a(B_+(t_1\ldots t_\mu))&=0,\quad \mu\geq 2 \\
 \mathbf{F}_a(B_+(t))&=\mathbf{a}(t),\ \forall t\in T_O
 \end{align*}
\end{lemma}

\begin{lemma}
Let $G\in V$ be any vector field with B-series of the form
$$
     hG = \sum_{t\in T_O} h^{|t|-1}\mathbf{G}(t) \mathbf{F}_x(B_+(t)),\qquad \mathbf{G}(\ab)=0.
$$
Then its $h$-flow $\exp(hG)$ has again a B-series $B(\mathbf{g},x)$ where
\begin{align*}
g(\emptyset) &= 1, \\
\mathbf{g}(B_+(t_1\ldots t_\mu))&=\frac1{\mu!}\mathbf{G}(t_1)\cdots\mathbf{G}(t_\mu)
\end{align*}
\end{lemma}
In order to obtain a systematic development of the order conditions of the commutator-free schemes of 
Section~\ref{subsubsec:commfree}, consider the definition \eqref{commfree1}-\eqref{commfree2}, and define
$$
Y_{r,0} = x\in M\quad\mbox{and}\quad Y_{r,j} = \exp\left(\sum_k \alpha_{r,j}^k F_k\right) Y_{r,j-1},\quad j=1,\ldots, J_r,
$$
in this way $Y_{r,J_r}=Y_r$ in \eqref{commfree1}. For the occasion, we unify the notation by setting
$\alpha_{s+1,j}^k := \beta_j^k$ for $1\leq k\leq s$ and $1\leq j\leq J$. Then we may also write $Y_{s+1}$ for $y_1$ in \eqref{commfree2} and we can write down the order conditions for the commutator free Lie group methods.
\begin{theorem}\label{theo:ocf}
With the above definitions,  the quantities $Y_{r,j},\ 1\leq r\leq s+1,\ 1\leq j\leq J_r$ all have B-series $B(\mathbf{Y}_{r,j},x)$ defined through the following formulas
\begin{align*}
\mathbf{Y}_{r,j}(\ab)&=1, \\
\mathbf{Y}_{r,0}(t)&=0,\quad \forall t: |t|>1,\\[1mm]
\mathbf{Y}_{r,j}(B_+(t_1\ldots t_\mu))&=\sum_{k=0}^\mu \mathbf{Y}_{r,j-1}(B_+(t_1\ldots t_k))\cdot\mathbf{b}_{r,j}(B_+(t_{k+1}\ldots t_\mu)),\\
\mathbf{b}_{r,j}(\ab)&=1,\quad 1\leq r\leq s+1,\ 1\leq j\leq J_r,\\
\mathbf{b}_{r,j}(B_+(t_1\ldots t_\mu))&=\frac{1}{\mu!}\mathbf{G}_{r,j}(t_1)\cdots\mathbf{G}_{r,j}(t_\mu),\\
\mathbf{G}_{r,j}(t) &=\sum_{k=1}^s \alpha_{r,j}^k\mathbf{Y}_{k,J_r}(t)
\end{align*}
\end{theorem}
It now have a B-series for the numerical solution and the B-series of the exact solution is given by \eqref{Bseries:exact} so it follows that
\begin{corollary}
A commutator-free Lie group method has order of consistency $q$ if and only if
$$
   \mathbf{Y}_{s+1,J}(t) = \frac{\alpha(t)}{(|t|-1)!},\quad \forall t: |t|\leq q+1.
$$
Here $\mathbf{Y}_{s+1,J}$ is found from Theorem~\ref{theo:ocf} and $\alpha(t)$ is given in Proposition~\ref{prop:alpha}.
\end{corollary}

A usual strategy for deriving methods of  a given order of consistency is to first consider the classical order conditions for Runge--Kutta methods since these of course must be satisfied for the classical Runge-Kutta coefficients defined in \eqref{classcoeff}. Since the computation of flows (exponential) is usually the most expensive operation, one next seeks the smallest possible number of exponentials per step, i.e. let each $J_r$ be as possible while leaving enough free parameters to solve the remaining non-classical conditions $\alpha_{r,j}^k$ and $\beta_j^k$.

\subsection{Selecting a minimal set of conditions.}
The conditions arising from each ordered rooted tree are not independent.
The B-series \eqref{Bseries} we consider here are representing objects of different kinds, such as maps and vector fields. 
This means that the series representing these objects are subsets of $\mathcal{T}_O$ and we briefly characterise these subsets here in a purely algebraic fashion. Then we can better understand how to select a minimal set of conditions to solve.

Let $A_- = \{t\in T_O:t=B_+(u), u\in T_O\}$ and set $A=A_-\cup\{\ab\}$. Any tree $t=B_+(t_1\ldots t_\mu)$ can be considered as a word of the alphabet $A$ in the sense that it is formed by the finite sequence $B_+(t_1),\ldots, B_+(t_\mu)$. With the concatenation product \eqref{concat} we get all of $T_O$ as the free monoid over the set $A$ with identity element $\ab$. This structure is then extended to $\mathcal{T}_O$ as an associative $\mathbb{R}$-algebra. The elements of $\mathcal{T}_O$ are the formal series on $T_O$ and we denote by $(P,t)\in\mathbb{R}$ the coefficient of the tree $t$ in the series $P$. The product of two series $S$ and $T$ is defined to be the series with coefficient
$$
(ST,t) = \sum_{t=u\cdot v} (S,u)(T,v).
$$
Next, notice that the commutator-free methods considered here are derived by composing flows of linear combinations of frozen vector fields. From Lemma~\ref{lemma:BseriesFrozen} we see that that these linear combinations (scaled by $h$) have expansions whose coefficents vanish on trees not belonging to $A_-$. On tfhe other hand, the composition of exponentials can be written as a single exponential of a series via the Baker--Campbell--Hausdorff formula, and the resulting series belongs to the free Lie algebra on the set $A_-\subset T_O$.

There are three important subsets of $\mathcal{T}_O$
\begin{itemize}
\item the subspace $\g\in\mathcal{T}_O$ which is the free Lie algebra on the set $A$
\item $V\subset\g$ is the subspace of $\g$ consisting of series $S$ such that
$$
     (S,t) = 0\quad\mbox{whenever}\ t\not\in A_-   
$$
\item $G$ is the group of formal exponential series $T=\exp(S),\; S\in\g$. These series are such that $(T,\ab)=1$.
\end{itemize}
The combinatorial properties of ordered rooted trees and free Lie algebras are by now well understood,  and many results hinge on on the Poincare-Birkhoff-Witt theorem, see e.g. \cite[Ch. 5]{humphreys72itl}. The space $\mathcal{T}_O$ has  a natural grading arising from the number of nodes in each ordered rooted tree, we can define
$$
      \mathcal{T}_O^q = \mbox{span}\{t: |t|=q+1\}\quad\mbox{thus}\ \mathcal{T}_O=\coprod_{n\geq 0}\mathcal{T}_O^n
$$
and similarly
$$
     \g = \coprod_{n\geq 0}\g^n,\quad \g^n = \g\cup\mathcal{T}_O^n
$$
The dimension of $\mathcal{T}_O^n$ is the Catalan number 
$$
\dim \mathcal{T}_O^n = \frac1{n+1}\left(\begin{array}{c} 2n\\ n\end{array}\right)
$$
The next result is well-known its proof can be found for instance in \cite{owren06ocf}
\begin{theorem}\label{theo:dimgn}
\[
    \dim\g_n = \nu_n=\frac{1}{2n}\sum_{d|n}\mu(d)\binom{2n/d}{n/d}
\]
where $\mu(d)$ is the M\"{o}bius function defined for any positive
integer as $\mu(1)=1$, $\mu(d)=(-1)^p$ when $d$ is the product of $p$
distinct primes, and $\mu(d)=0$ otherwise. The sum is over all
positive integers which divide $n$, including 1 and $n$.
\end{theorem}
We present a table over the numbers $\mathcal{T}_O^n$, $\g^n$ and $c^n$, the last one being the dimensions of the graded components of the unordered trees, that count the number of order conditions for classical RK methods
$$
\begin{array}{|c|xxxxxxx|}
\hline 
  n     &  1 & 2 & 3 & 4 & 5 & 6 & 7 \\ \hline
\mathcal{T}_O^n & 1 & 2 & 5 & 14 & 42 & 132 & 429 \\ \hline
\g^n &  1 & 1 & 3 & 8 & 25 & 75 & 245 \\ \hline
c^n  & 1 & 1 & 2 & 4 & 9 & 20 & 48  \\ \hline
\end{array}
$$
So in fact, the numbers $\g^n$ gives the number of order conditions to be considered for each order for commutator-free methods and a possible strategy would be to pick $\g^n$ independent conditions out of the $\mathcal{T}_O^n$ found from 
Theorem~\ref{theo:ocf}. It should be observed that the dependency between conditions corresponding to ordered rooted trees arise amongst trees that share the same (unordered) set of subtrees, such as $\aaabbabb$ and $\aabaabbb$, in fact the condition corresponding to precisely one of these two trees can be discarded given that conditions of lower order are included.
Using classical theory of free Lie algebras, one may characterize this dependency by a generalized Witt formula counting, for a given tree t, the dimension of the subspace of $\g$ spanned by the set of trees obtained from permuting the subtrees of $t$. Consider the equivalence class $[t]$ characterized by a set of $\nu$ distinct subtrees, $t_i\in T_O,\; i=1,\ldots,\nu$, where there are exactly $\kappa_i$ occurrences of $t_i$. The dimension of the subspace spanned by trees in $[t]$ can be derived from a formula in Bourbaki \cite{bourbaki75lga}.
$$
      c(\kappa) = \frac{1}{|\kappa|} \sum_{d|\kappa} \mu(d)\frac{(|\kappa|/d)!}{(\kappa/d)!},
$$
in other words, the dimensions depends only on the number of copies of each subtree and not on the subtrees themselves.
For convenience, we give a few examples
\begin{equation}\label{eq10}%
\begin{array}{lcll}
   c(n)   & = & 0, & n>1 \\
   c(n,1) &=& 1,\quad& n>0,\\
   c(n,1,1)&=& n+1,& n>0,\\
   c(n,2)&=& \lfloor \frac{n+1}{2}\rfloor & n>0.\\
   c(1,\ldots,1) &=& (r-1)!,\quad & (\mbox{$r$ distinct subtrees})
\end{array}
\end{equation}
A detailed analysis of how methods of orders up to four are constructed can be found in \cite{owren06ocf}. Schemes of order 4 with 4 stages can be constructed with two exponentials in the fourth stage and the update stage as in \eqref{cfexample}.
In this example, there is also a clever reuse of an exponential such that the total number of flow calculation is effectively 5.
Yet another result from \cite{owren06ocf} shows that in the case of two exponentials, i.e. $J_r=2$ the coefficients of stage $r$ is only involved linearly in the order conditions.

As far as we know, no explicit commutator-free method of order five or higher has been derived at this point. A complication is then that one needs to have stages (including the final update) with at least three exponentials, and no simplification of the order conditions similar to the $J_r=2$-case has been found.

\section{ Symplectic Lie group integrators} \label{sec:symplectic}
It is not clear whether there exist Lie group integrators of Munthe-Kaas  or commutator-free type which are symplectic for an arbitrary symplectic manifold. Recently, McLachlan et al. \cite{mclachlan14amc} found an elegant adaptation to the classical midpoint rule to make it a symplectic integrator on product of 2-spheres. It is however relatively easy to find symplectic integrators on cotangent bundles of manifolds, an by looking at the special case where $M=T^*G$ for a Lie group $G$ one can obtain symplectic integrators which are rather similar in form to partitioned Munthe-Kaas methods as defined by Eng{\o} in \cite{engo03prk}. The approach we described here was introduced by Celledoni et al. \cite{celledoni14ait} and by using ideas from Bou-Rabee and Marsden \cite{bou-rabee09hpi}, Bogfjellmo and Marthinsen \cite{bogfjellmo15hos} extended it to high order partitioned symplectic Lie group integrators both in the Munthe-Kaas and Crouch-Grossman formats. We briefly describe the setting for these symplectic integrators.

The first step is to replace $T^*G$ by $\mathbf{G}:=G\ltimes\g^*$ via right trivialisation, meaning that any $p_g\in T_g^*G$ is represented 
as the tuple ($g,\mu)$ where
$\mu=R_g^*p_g$. We use the notation $R_{g*}v=TR_gv$ for any $v\in TG$ and similarly $R_g^*$ for the adjoint operator
such that $\langle R_g^*p,v\rangle=\langle p,R_{g*}v\rangle$ for any $p\in T^*G$, $v\in TG$.
Next, we lift the group structure from $G$ to $\mathbf{G}$ through
$$
(g_1,\mu_1)\cdot(g_2,\mu_2) = (g_1\cdot g_2, \mu_1+\Ad_{g_1^{-1}}^*\mu_2), \qquad  1_{\mathbf G} =(1_G,0_{\g^*}),
$$
where $\Ad_g^*$ is defined in \eqref{eq:coadjoint-action}. Similarly, the tangent map of right multiplication extends as
\[
TR_{(g,\mu)}(R_{h*}\,\zeta, \nu) = (R_{hg*}\;\zeta, \nu-\ad_\zeta^*\,\Ad_{h^{-1}}^*\mu),\quad
g,h \in G,\ \zeta\in\g,\ \mu,\nu\in\g^*.
\]
Of particular interest is the restriction of $TR_{(g,\mu)}$ to $T_1\mathbf{G}\cong \g \times \g^*$,
\[
    T_1R_{(g,\mu)}(\zeta,\nu) = (R_{g*}\zeta, \nu - \ad_\zeta^*\mu). \]
The natural symplectic form on
$T^*G$  is defined as
\[
\Omega_{(g,p_g)}((\delta v_1, \delta \pi_1),(\delta v_2,\delta\pi_2))
=\langle \delta \pi_2, \delta v_1\rangle - \langle \delta \pi_1, \delta v_2\rangle,
\]
and by right trivialisation it may be pulled back to $\mathbf{G}$ and then takes the form
\begin{equation} \label{eq:sympform}
      \omega_{(g,\mu)}( (R_{g*} \xi_1, \delta\nu_1), (R_{g*}\xi_2, \delta\nu_2))
      = \langle\delta\nu_2,\xi_1 \rangle - \langle\delta\nu_1, \xi_2 \rangle - \langle\mu, [\xi_1,\xi_2]\rangle,
\end{equation}
where  $\xi_1,\xi_2 \in \g$.
The presentation of differential equations on $T^*G$ as in \eqref{F_and_rho} is now achieved via the action by left multiplication, meaning that any vector field $F\in\mathcal{X}(\mathbf{G})$ is expressed by means of a map $f\colon \mathbf{G} \rightarrow T_1 \mathbf{G}$,
\begin{equation} \label{eq:Fpres}
        F(g,\mu) = T_1 R_{(g,\mu)} f(g,\mu) = (R_{g*} f_1, f_2-\ad_{f_1}^*\mu),
\end{equation}
where $f_1=f_1(g,\mu)\in\g$, $f_2=f_2(g,\mu)\in\g^*$ are the two components of $f$.
We are particularly interested in the case that $F$ is a Hamiltonian vector field
which means that $F$ satisfies the relation
\begin{equation} \label{eq:iF}
    F\intder \omega = \d H,
\end{equation}
for some Hamiltonian function $H\colon T^*G\rightarrow\R$ and $\intder$ is the interior product defined as $F \intder \omega (X) := \omega(F, X)$ for any vector field $X$.
From now on we let $H \colon \mathbf{G} \to \R$ denote the trivialised Hamiltonian.
A simple calculation using 
\eqref{eq:sympform}, \eqref{eq:Fpres} and \eqref{eq:iF} shows that the corresponding map $f$ for such a Hamiltonian vector field is
\begin{equation} \label{rkmk-f}
     f(g,\mu) = \left(\frac{\partial H}{\partial\mu}(g,\mu), -R_g^*\frac{\partial H}{\partial g}(g,\mu)\right).
\end{equation}

\subsection{Variational integrators on Lie groups}

A popular way to derive symplectic integrators in general is through the discretisation of a variational principle, see e.g. Marsden and West \cite{marsden01dma}. The exact solution to the mechanical system is the function $g(t)$ which minimises the action
\begin{equation} \label{actionintegral}
      S[g] = \int_{t_0}^{t_1} L(g,\dot{g})\,\d t.
\end{equation}
The function 
$L:TM\rightarrow \mathbb{R}$ is called the Lagrangian.
Taking the variation of this integral yields the Euler--Lagrange equations
$$
  \frac{\d}{\d t}   \frac{\partial L}{\partial \dot{g}} = \frac{\partial L}{\partial g}.
$$
It is usually more conventient to trivialise $L$ by defining $\ell(g,\xi)=L(g,\R_{g*}\xi)$, $\ell: G\ltimes\g\rightarrow\mathbb{R}$. In this case it is an easy exercise to derive the corresponding version of the Euler--Lagrange equations as
$$
\frac{\d}{\d t} \frac{\partial\ell}{\partial\xi} = R_g^*\frac{\partial\ell}{\partial g} - \ad_{\xi}^*\frac{\partial\ell}{\partial\xi}.
$$
At this point it is customary to introduce the Legendre transformation, defining 
\begin{equation} \label{legendretransformation}
\mu=\frac{\partial\ell}{\partial\xi}(g,\xi)\in\g^*\quad\Rightarrow\quad \xi=\iota(g,\mu)\ \text{where}\ \iota: G\times\g^*\rightarrow\g,
\end{equation}
assuming that the left equation can be solved for $\xi$. The corresponding differential equations for $g$ and $\mu$ would be
$$
      \dot{g} = R_{g*}\iota(g,\mu),\qquad \dot{\mu} = R_g^*\frac{\partial\ell}{\partial g}(g,\iota(g,\mu))
      -\ad_{\iota(g,\mu)}^* \mu.
$$
which agrees with the formulation \eqref{eq:Fpres} with 
\begin{equation}\label{f-explicit-form}
f_1(g,\mu)=\iota(g,\mu)\quad\text{and}\quad
f_2(g,\mu)=R_g^* \frac{\partial\ell}{\partial g}(g,\iota(g,\mu)).
\end{equation}
Variational integrators are derived by extremising an approximation to \eqref{actionintegral},
$$
   S_d\left(\{g_k\}_{k=0}^{N-1}\right) = h\sum_{k=0}^{N-1} L_d(g_k,g_{k+1}),
$$
with respect to the discrete trajectory of points $g_k\approx g(t_k)$. The function $L_d(x,y)$ is called \emph{the discrete Lagrangian.}
We compute the variation
\begin{multline*}
\delta S_d=h\sum_{k=1}^{N-1}
\langle
D_1L_d(g_k,g_{k+1})+D_2L_d(g_{k-1},g_k), \delta g_k
\rangle
\\+
h\langle D_1L_d(g_0,g_1), \delta g_0 \rangle
+h \langle D_2L_d(g_{N-1}, g_N),\delta g_N\rangle.
 \end{multline*}
Leaving the end points fixed amounts to setting $\delta g_0=\delta g_N=0$ and this leads to the discrete version of the Euler-Lagrange equations (DEL)
$$
      D_1 L_d(g_k,g_{k+1}) + D_2 L_d(g_{k-1},g_k)=0,\quad k=1,\ldots,N-1.
$$
We now present an example of how $L_d(g_k,g_{k+1})$ can be defined. Let $\xi^{-1}: \g\rightarrow G$ be a map which satisfies 
\begin{enumerate}
\item $\xi^{-1}$ is a local diffeomorphism
\item $\xi^{-1}(0) = 1$
\item $T_0\xi^{-1} = \mathrm{Id}_{\g}$
\end{enumerate}
For any $\eta\in\g$  define the curve $g(t)=\xi^{-1}(t\eta) g_k$ on the Lie group. Clearly, we get $g(0)=g_k$ and
by choosing
$$
     \eta = \xi(g_{k+1}g_k^{-1}) = \xi(\Delta_k),\quad \Delta_k := g_{k+1} g_k^{-1},
$$
we get  $g(1)=g_{k+1}$.  As for the second argument of the continuous Lagrangian, we  compute
 $\dot{g}(t)=R_{g_k*}\circ T_{t\eta}\xi^{-1}\circ\eta$. Taking $\dot{g}(0)$ as an approximation to the argument $\dot{g}$ of $L$, we find by 3.\ the approximation 
$\dot{g}(0)=R_{g_k*}\xi(g_{k+1}g_k^{-1})=R_{g_k*}\xi(\Delta_k):=R_{g_k*}\xi_k$.
So a possible approximation $L_d(g_k,g_{k+1})$ could be
$$
     L_d(g_k,g_{k+1}) = L(g_k, R_{g_k*}\xi_k) := \ell(g_k,\xi_k).
$$
We can compute the variation of the action sum
\begin{multline*}
\delta S_d = h\sum_{k=0}^{N-1}\bigg(
\langle \frac{\partial\ell}{\partial g}(g_k,\xi_k)
-R_{g_k^{-1}}^*\Ad_{\Delta_k}^*\mathrm{d}\xi_{\Delta_k}^*\frac{\partial\ell}{\partial\xi}(g_k,\xi_k),\delta g_k\rangle \\
+\langle
R_{g_{k+1}^{-1}}^*\mathrm{d}\xi_{\Delta_k}^*\frac{\partial\ell}{\partial\xi}(g_k,\xi_k),
\delta g_{k+1}
\rangle
\bigg).
\end{multline*}
The DEL equations are thus
$$
     R_{g_k}^*\frac{\partial\ell}{\partial g}(g_k,\xi_k)
     -\Ad_{\Delta_k}^*\mathrm{d}\xi_{\Delta_k}^*\frac{\partial\ell}{\partial \xi}(g_k,\xi_k)
     +\mathrm{d}\xi_{\Delta_{k-1}}^* \frac{\partial\ell}{\partial\xi}(g_{k-1},\xi_{k-1})=0,\quad
     1\leq k\leq N-1.
$$
The trivialised discrete Legendre transformations can be seen as maps 
$\mathbb{F}L_d^{\pm}: G\times G\rightarrow G\times\g^*$.
\begin{align}
\mathbb{F}L_d^+(g,g') &=(g',R_{g'}^*D_2L_d(g,g'))=(g',
\mathrm{d}\xi_\Delta^*\frac{\partial\ell}{\partial\xi}(g,\xi(\Delta))) \label{eq:FLplus}, \\
\mathbb{F}L_d^-(g,g') &=(g,-R_{g}^*D_1L_d(g,g'))=
(g,-R_{g}^*\frac{\partial\ell}{\partial g}(g,\xi(\Delta))
+\Ad_\Delta^*\mathrm{d}\xi_{\Delta}^*\frac{\partial\ell}{\partial\xi}(g,\xi(\Delta)))\label{eq:FLminus}, \\
\Delta &= g'g^{-1}.
\end{align}
One way of interpreting a method on the trivialised contangent bundle $G\rtimes\g^*$ is
the following
\begin{enumerate}
\item
Assume that $(g_k,\mu_k)$ is known. Then compute $g_{k+1}$ by solving the equation
$$
       \mathbb{F}L_d^-(g_k,g_{k+1}) = \mu_k.
$$
\item 
Next solve for $\mu_{k+1}$ explicitly by
$$
      \mu_{k+1} = \mathbb{F}L_d^+(g_k,g_{k+1}).
$$
\end{enumerate}
What we would like is to replace the occurrences of the Lagrangian $\ell(g,\xi)$ by
the functions used in the RKMK formulation \eqref{eq:Fpres}, \eqref{rkmk-f}.
A plausible start is to define the variable $\bar{\mu}_k$ by
$$
     \bar{\mu}_k := \frac{\partial\ell}{\partial\xi}(g_k,\xi(g_{k+1}g_k^{-1})).
$$
The continuous Legendre transformation \eqref{legendretransformation}  yields
$$
       \xi(g_{k+1}g_k^{-1})=\xi(\Delta_k) = \iota(g_k,\bar{\mu}_k)=f_1(g_k,\bar{\mu}_k),
$$
thus
$$
       g_{k+1} = \xi^{-1}(\iota(g_k,\bar{\mu}_k))\cdot g_k.
$$
To find an equation for $\bar{\mu}_k$ we need to consider $\mathbb{F}L_d^-(g_k,g_{k+1})$ as described in point 1.\ above.
Notice, again from \eqref{f-explicit-form} that 
$$
     R_{g_k}^*\frac{\partial\ell}{\partial g}(g_k,\xi(\Delta_k)) = f_2(g_k,\bar{\mu}_k).
$$
From \eqref{eq:FLminus} with $g=g_k$ and $g'=g_{k+1}$ we get
$$
    \mu_k = -f_2(g_k,\bar{\mu}_k) + \Ad_{\Delta_k}^*\mathrm{d}\xi_{\Delta_k}^*\bar{\mu}_k.
$$
Finally, we need only use \eqref{eq:FLplus} to get $\mu_{k+1}$. All equations for one step can be summarized as
\begin{align*}
    g_{k+1} &= \xi^{-1}(f_1(g_k,\bar{\mu}_k))\cdot g_k, \\
     \mu_k & = -f_2(g_k,\bar{\mu}_k) + \Ad_{\Delta_k}^*\mathrm{d}\xi_{\Delta_k}^*\bar{\mu}_k,  \\
     \mu_{k+1} &= \mathrm{d}\xi_{\Delta_k}^* \bar{\mu}_k,
\end{align*}
where as before $\Delta_k=g_{k+1}g_k^{-1}$. We could use the  map
$\tau=\xi^{-1}$ instead, recalling that for $u=\xi(g)$ one has $\mathrm{d}\xi_g=(\mathrm{d}\xi^{-1}_u)^{-1}$. One way of writing the resulting method would be
\begin{equation} \label{eq:eulerBO}
\begin{split}
    g_{k+1} &= \tau(f_1(g_k,\mathrm{d}\tau_{\xi_k}^*\mu_{k+1}))\cdot g_k, \\
    \mu_{k+1} &= \Ad_{\Delta_k^{-1}}^*(\mu_k+f_2(g_k,\mathrm{d}\tau_{\xi_k}^*\mu_{k+1})),
\end{split}
\end{equation}
where $\xi_k=\xi(\Delta_k)$, i.e. $\Delta_k=\tau(\xi_k)$.
It is already well-known that  a scheme derived from such a variational principle leads to a symplectic method, see e.g. Marsden and West \cite{marsden01dma}. By replacing the discrete Lagrangian and action sum by other more advanced approximations, one can obtain various different variants of symplectic integrators on Lie groups, see e.g. \cite{bogfjellmo15hos,celledoni14ait}.

\section{ Preservation of first integrals} \label{sec:intpres}
There has been a significant interest over the last decades in constructing integrators which preserve one or more first integrals, such as energy or momentum. The reader who is interested in this topic should consult the pioneering paper by Gonzalez \cite{gonzalez96tia}, but also McLachlan et al. \cite{mclachlan99giu} and the more recent work \cite{quispel08anc} and for preserving multiple first integrals simultaneously, see \cite{dahlby11pmf}. A key concept in integral preserving methods is that of discrete gradients, and in \cite{celledoni14pfi} these concepts were extended to Lie groups and retraction manifolds.

We shall begin by considering the case of a Lie group $G$ and define a first integral to be any differentiable function $H: G\rightarrow\mathbb{R}$ which is invariant on solutions
$$
  \frac{\d}{{\d}t} H(y(t))   = \langle{\d}H,F\rangle = 0,
$$
where we have introduced a duality pairing $\langle\cdot,\cdot\rangle$ between vector fields and one-forms.
To any differential equation on a Lie group having a first integral $H$, there exists a bivector (dual two-form) $\omega$
such that
\begin{equation} \label{bivector_ode}
\dot{y} = F(y) = \omega(\d H,\cdot) = \d H \intder \omega.
\end{equation}
An explicit formula for $\omega$ can be given in the case when $M$ is a Riemannian manifold. The gradient vector field is defined at the point $x$ through the relation $\langle \d H,v\rangle_x=(\grad H|_x,v)_x$ for every $v\in T_xM$ where $(\cdot,\cdot)$ is the Riemannian inner product. An example of a bivector to be used in \eqref{bivector_ode} is then given by
$$
   \omega = \frac{\grad H\wedge F}{\|\grad H\|^2}.
$$
One should note that the bivector $\omega$ used to express the differential equation is not unique for a given ODE vector field $F$, but the choice of bivector will affect the resulting numerical method. The formulation \eqref{bivector_ode} can easily be generalised to the case of $k$ invariants, $H_1,\ldots, H_k$. In this case we replace the bivector by a $(k+1)$-vector and write the equation as
$$
   \dot{y} = F(y) = \omega(\d H_1,\ldots,\d H_k,\cdot),
$$
and again, for Riemannian manifolds, we can define $\omega$ as
$$
    \omega = \frac{\omega_0\wedge F}{\omega_0(\d H_1,\ldots,\d H_k)}\quad\mbox{where}\quad
    \omega_0 = \grad H_1 \wedge\cdots\wedge\grad{H}_k.
$$

\paragraph{\bf Integral preserving schemes on Lie groups}
Let $G$ be a Lie group with Lie algebra $\g$, and define for each $g\in G$, the right multiplication operator $R_g:G\rightarrow G$ by $R_g(h)=h\cdot g$.
\begin{definition}
Let $H\in\mathcal{F}(G)$. We define the trivialised discrete differential $\bar{\d}H: G\times G \rightarrow \g^*$  as any map that satisfies the conditions
\begin{align*}
H(v)-H(u) &= \langle\bar{\d}H(u,v),\log(v\cdot u^{-1})\rangle, \\
\bar{\d}H(g,g) &= R_g^*\d H_g,\qquad \forall g\in G.
\end{align*}
\end{definition}
We also need a trivialised approximation to the bivector $\omega$ in \eqref{bivector_ode}. For every pair of points
$(u,v)\in G\times G$, we define an exterior 2-form on the linear space $\g^*$,
$\bar{\omega}:G\times G\rightarrow \Lambda^2(\g^*)$, satisfying the consistency condition
 \begin{equation*}
     \bar{\omega}(g,g)(R_g^*\alpha, R_g^*\beta) = \omega_g(\alpha,\beta),\qquad \forall g\in G,\quad\forall \alpha,\beta\in T_g^*G.
 \end{equation*}
 For practical purposes, $\bar\omega$ needs only to be defined in some suitable neighborhood of the diagonal subset
 $\{(g,g):g\in G\}$. We can now write the numerical integral preserving method as
 $$
 g^{i+1} = \exp(h\zeta(g^i,g^{i+1}))\cdot g^i,\qquad \zeta(g^i,g^{i+1})=\bar\d H(g^i,g^{i+1})\intder\bar\omega(g^i,g^{i+1}).
 $$
That the integral $H$ is exactly preserved is seen through the simple calculation
\begin{equation*}
 \begin{split}
     H(g^{i+1})-H(g^i)&=\langle \bar\d H(g^i,g^{i+1}), \log(g^{i+1}(g^i)^{-1})\rangle
     =h\langle \bar\d H(g^i,g^{i+1}), \zeta(g^i,g^{i+1})\rangle \\
     &= \bar{\omega}(g^i,g^{i+1})(\bar\d H(g^i,g^{i+1}),\bar\d H(g^i,g^{i+1}))=0.
\end{split}
\end{equation*}

\paragraph{\bf Examples of trivialised discrete differentials}

The first example has a counterpart on Euclidean space sometimes referred to as the \emph{Averaged Vector Field} (AVF) gradient. It is defined as
\begin{equation} \label{LG-avf}
\bar\d H(u,v) = \int_0^1 R_{\ell(\xi)}^*\d H_{\ell(\xi)}\,\d\xi,\quad
\ell(\xi)=\exp(\xi\log(v\cdot u^{-1})).
\end{equation}
Note that $\bar\d H(u,v)=\bar\d H(v,u)$. This trivialised discrete differential has the disadvantage that it can rarely be computed exactly a priori for general groups, although when $G$ is Euclidean space it reduces to the standard AVF discrete gradient which has a wide range of applications.

Gonzalez \cite{gonzalez96tia} introduced a discrete gradient for Euclidean spaces which is oftentimes referred to as the midpoint discrete gradient. In the setting we use here, we need to introduce an inner product $(\cdot,\cdot)$ on the Lie algebra to generalise it to arbitrary Lie groups.
We apply ``index lowering" to any element $\eta\in\g$ by defining $\eta^\flat\in\g^*$ to be  the unique element satisfying
$\langle\eta^\flat, \zeta\rangle=(\eta,\zeta)$ for all $\zeta\in\g$.  We let
\begin{equation} \label{LG-gmp}
\bar\d H(u,v) = R_c^*\d H|_c + \frac{H(v)-H(u)-\langle R_c^* \d H|_c,\eta\rangle}{(\eta,\eta)}\,\eta^\flat,\quad
\eta = \log(v\cdot u^{-1}),
\end{equation}
where $c\in G$, is some point typically near $u$ and $v$. One may for instance choose $c=\exp(\eta/2)\cdot u$, which implies symmetry, i.e. $\bar\d H(u,v)=\bar\d H(v,u)$.

\paragraph{\bf Integral preserving schemes on a manifold with a retraction}

What we present here is a basic and straightforward approach introduced in \cite{celledoni14pfi}, but clearly there are other strategies that can be used. We use retractions as introduced on page~\pageref{retractions}. Recall that the retraction restricted to the tangent space $T_xM$ is denoted $\phi_x$ and is a diffeomorphism from some neighborhood $\mathcal{U}_x$ of $0_x\in T_xM$ into a subset 
$\mathcal{W}_x$ of of $M$ containing $x$. We tacitly assume the necessary restrictions on the integration  step size $h$ to ensure that both the initial and terminal points are contained in $\mathcal{W}_x$ for each time step. We also assume to have at our disposal a map $c$ defined on some open subset of $M\times M$ containing all diagonal points $(x,x)$, for which $c(x,y)\in M$.
The differential equation is written in terms of a bivector $\omega$ and a first integral $H$ as in \eqref{bivector_ode}. We introduce an approximate bivector $\bar\omega(x,y)$ such that
$$
\bar\omega(x,x)(v,w) = \left.\omega\right|_x,\quad\forall x\in M.
$$
Contrary to the Lie group case we no longer assume a global trivialisation, so we introduce the \emph{discrete differential} of a function $H$. To any pair of points $(x,y)\in M\times M$ we associate the covector $\bar\d H(x,y)\in T_{c(x,y)}^*M$ satisfying the relations
\begin{align*}
    H(y)-H(x) &= \langle\d H(x,y),\phi_c^{-1}(y)-\phi_c^{-1}(x)\rangle, \\[1mm]
    \bar\d H(x,x) &= \left.\d H\right|_x,\quad\forall x\in M.
\end{align*}
where $c=c(x,y)$ is the map introduced above. The integrator on $M$ is now defined as
$$
  y^{n+1} = \phi_c(W(y^n,y^{n+1})),\quad W(y^n, y^{n+1})=\phi_c^{-1}(y^n)+h\,\bar\d H(y^n,y^{n+1})\intder\bar\omega(y^n,y^{n+1}).
$$
One can easily see that this method is symmetric
 if the following three conditions are satisfied:
\begin{enumerate}
\item The map $c$ is symmetric, i.e. $c(x,y)=c(y,x)$ for all $x$ and $y$.  \label{csym}
\item The discrete differential is symmetric in the sense that $\bar\d H(x,y)=\bar\d H(y,x)$. \label{dbsym}
\item The bivector $\bar{\omega}$ is symmetric in $x$ and $y$: $\bar{\omega}(x,y)=\bar\omega(y,x)$. \label{omsym}
\end{enumerate}
The condition \ref{csym}) can be achieved by solving the equation
\begin{equation}\label{eq:symccond}
     \phi_c^{-1}(x) + \phi_c^{-1}(y) = 0,
\end{equation}
with respect to $c$.

Both of the trivialised discrete differentials \eqref{LG-avf} and \eqref{LG-gmp} have corresponding versions with retractions, in the former case, we write $\gamma_\xi=(1-\xi)v+\xi w$ where
$x=\phi_c(v)$, $y=\phi_c(w)$. Then
\begin{equation} \label{eq:avfgenman}
\bar\d H(x,y) = \int_0^1\phi_c^* \left.\d H\right|_{\phi_c(\gamma_\xi)} \, \d\xi.
\end{equation}
Similarly, assuming that $M$ is Riemannian, we can define the following counterpart to the Gonzalez midpoint discrete gradient
\begin{equation} \label{eq:gonzgenman}
\bar\d H(x,y) = \d H|_c + \frac{H(y)-H(x)-\langle \d H|_c,\eta\rangle}{(\eta,\eta)_c}\,\eta^\flat,\quad
 \eta = \phi_c^{-1}(y)-\phi_c^{-1}(x) \in T_cM.
\end{equation}
where we may require that $c(x,y)$ satisfies \eqref{eq:symccond} for the method to be symmetric.
For clarity, we include an example taken from \cite{celledoni14pfi}

\begin{example}
We consider the sphere $M=S^{n-1}$ where we represent its points as vectors in $\mathbb{R}^n$ of unit length,
$\|x\|_2=1$. The tangent space at $x$ is then identified with the set of vectors in $\mathbb{R}^n$ orthogonal to $x$ with respect to the Euclidean inner product $(\cdot,\cdot)$.
 A retraction is
\begin{equation} \label{eq:rets2}
     \phi_x(v_x) = \frac{x+v_x}{\|x+v_x\|},
\end{equation}
its inverse is defined in the cone $\{y: (x,y)>0\}$ where
$$
     \phi_x^{-1}(y) = \frac{y}{( x, y)} - x.
$$
A symmetric map $c(x,y)$ satisfying \eqref{eq:symccond} is simply
\begin{equation} \label{eq:centerpoint}
      c(x,y) = \frac{x+y}{\|x+y\|_2},
\end{equation}
the geodesic midpoint between $x$ and $y$ in terms of the standard Riemannian metric on $S^{n-1}$.
We compute the tangent map of the retraction to be
$$
      T_u\phi_c = \frac{1}{\|c+u\|_2}\left(I-\frac{(c+u)\otimes(c+u)}{\|c+u\|_2^2}\right).
$$
As a toy problem, let us consider a mechanical system on $S^2$. Since the angular momentum in body coordinates for the free rigid body is of constant length, we may assume $(x,x)=1$ for all $x$ and we can model the problem as a dynamical system on the sphere. But in addition to this, the energy of the body i preserved,
i.e.
$$
      H(x) = \frac12(x, \mathbb{I}^{-1}x) = \frac12\left(\frac{x_1^2}{\mathbb{I}_1}+\frac{x_2^2}{\mathbb{I}_2}+\frac{x_3^2}{\mathbb{I}_3}\right),
$$
which we may take as the first integral to be preserved. Here the inertia tensor is $\mathbb{I}=\mbox{diag}(\mathbb{I}_1, \mathbb{I}_2, \mathbb{I}_3)$. The system of differential equations can be written as follows
\begin{align*}
\frac{\d x}{\d t} &=  \left.(\d H \intder \omega)\right|_x = x \times \mathbb{I}^{-1} x, \\
\omega|_x(\alpha,\beta) &= (x,\alpha\times\beta),
\end{align*}
where the righthand side in both equations refer to the representation in $\mathbb{R}^3$.
A symmetric consistent approximation to $\omega$ would be
$$
  \bar{\omega}(x,y)(\alpha,\beta) =(\frac{x+y}{2}, \alpha\times\beta).
$$
We write $\ell_\xi=c+\gamma_\xi$ with the notation in \eqref{eq:avfgenman}, this is a linear function of the scalar argument $\xi$.
and thus, $\phi_c(\gamma_\xi)=\ell_{\xi}/\|\ell_\xi\|$ from \eqref{eq:rets2}. We therefore derive for the AVF discrete gradient
 $$
\bar\d  H(x,y) = \int_0^1 \frac{1}{\|\ell_\xi\|}\left(\mathbb{I}^{-1}\phi_c(\gamma_\xi)  - (\phi_c(\gamma_\xi),\mathbb{I}^{-1}\phi_c(\gamma_\xi))
\phi_c(\gamma_\xi)\right)\;
\d\xi.
$$
This integral is somewhat complicated to solve analytically. Instead, we may consider the discrete gradient \eqref{eq:gonzgenman}
where we take as Riemannian metric the standard Euclidean inner product restricted to the tangent bundle of $S^2$.
We obtain the following version of the discrete differential in the chosen representation
$$
\bar\d  H(x,y)=  \frac{1}{\|m\|}\left(\mathbb{I}^{-1}m + \frac{\|m\|^2-1}{\|y-x\|^2}(H(y)-H(x))(y-x)\right),\quad m=\frac{x+y}{2}.
$$
The corresponding method is symmetric, thus of second order.
\end{example}

\bibliographystyle{plain}
\bibliography{brynbib,geom_int}

\end{document}